\documentclass[11pt]{elsarticle}

\usepackage{amsmath}
\numberwithin{equation}{section}
\usepackage{amssymb}
\usepackage{amsfonts,amscd,amsthm,amsbsy,bm,latexsym}
\usepackage{bookmark}
\usepackage{color,xcolor}
\usepackage{color}
\usepackage{graphicx}
\usepackage{subfigure}

\usepackage{amsthm}

\newtheorem{lemma}{Lemma}[section]
\newtheorem{theorem}{Theorem}[section]
\newtheorem{definition}{Definition}[section]
\newtheorem{proposition}{Proposition}[section]
\newtheorem{remark}{Remark}[section]
\newtheorem{corollary}{Corollary}[section]

\numberwithin{equation}{section}

\journal{ The Journal of Geometric Analysis}
\begin{document}

\begin{frontmatter}


\title{Well-posedness for the  semilinear wave equations with nonlinear damping on manifolds with conical degeneration}

\author[1]{Gongwei Liu \corref{cor1}}
\ead{gongweiliu@haut.edu.cn} 
\author[1]{Yi Peng}
\ead{pengyi0141@163.com}
\author[2]{Peng Li}
\ead{lipeng@haut.edu.cn}
\address[1]{School of Mathematics and  Statistics, Henan University of Technology, Zhengzhou 450001}
\address[2]{Institute for Complexity Science, Henan University of Technology, Zhengzhou 450001}

\begin{abstract}

This paper deals with a class of semilinear wave equation with nonlinear damping  term $|u_{t}|^{m-2}u_t $ and nonlinear source term $g(x)|u|^{p-2}u$ on the manifolds with conical singularities. Firstly, we prove the local existence and uniqueness of the  solution by the semigroup method. Secondly, we establish the global existence, the energy decay estimate and the blow-up in finite time of the solution with subcritial  ($E(0)<d$) and critial ($E(0)=d$) initial energy level by constructing potential wells. We also show that the solution is  global provided the damping dominates the source (that is $m\geq p$). Moreover, we prove the blow-up in finite time of the solution with arbitrary positive initial energy and give some necessary and sufficient condition for existing finite time blow-up solutions.
\end{abstract}


\end{frontmatter}


\section{Introduction}
\setcounter{equation}{0}
In this paper, we investigate the following semilinear wave equation with nonlinear damping and nonlinear source term on the manifolds with conical singularities
\begin{equation}\label{1.1}
\begin{cases}
u_{tt}-\Delta_{\mathbb{B}}u-\gamma V(x)u+|u_t|^{m-2}u_t=g(x)|u|^{p-2}u ,~~\,\,x\in int\mathbb{B}, \,t>0,\\[2mm]
 u(x, 0)=u_0(x),\,\,\, u_t(x, 0)=u_1(x),\,\, \,\,\,\,\,\,\qquad \qquad \qquad x\in int\mathbb{B},\\[2mm]
u(x,t)=0, ~~~~~~~~~~~~~~~~~~~~~\qquad \qquad \qquad \qquad \qquad x\in\partial\mathbb{B} ,\,t\geq0.
\end{cases}
\end{equation}
Here the domain  $\mathbb{B}=[0, 1)\times X$, $X$ is an $(n-1)$-dimensional compact smooth manifold ($n\geq3$), which is regarded as the local model near the conical points and $\partial \mathbb{B}=\{0\}\times X$. Near $\partial\mathbb{B}$, we use the local coordinates $x=(x_1, x')=(x_1,x_2,\cdots,x_n)$ for $x_1\in[0,1)$ and $x'\in X$.   Then the Fuchsian type Laplace operator is defined as $\Delta_{\mathbb{B}}:=(x_{1}\partial_{x_{1}})^{2} +\partial^{2}_{x_{2}}+\cdots+\partial^{2}_{x_{n}}$, which is an elliptic operator with conical degeneration on the boundary $x_1=0$, and the corresponding gradient operator is denoted by $\nabla_{\mathbb{B}}:=x_1\partial_{x_1}+\partial_{x_2}+\cdots+\partial_{x_n}$. $V(x)$ is a positive  potential function which can be unbounded on the cone manifold $\mathbb{B}$ and satisfies the cone Hardy's inequality (see the choice in Lemma 2.1). The initial datas $u_0(x)$ and $u_1(x)$ belong to the cone Sobolev space $\mathcal{H}_{2,0}^{1,\frac{n}{2}}(\mathbb{B})$ and  the cone Lebesgue space $L^{\frac{n}{2}}_2(\mathbb{B}) $, respectively. We refer to the references \cite{chen2012cone,chen2011existence,schulze1998boundary} for the detail research of manifold with conical singularities and the corresponding cone Sobolev spaces.

For the problem \eqref{1.1}, we always assume the following hypotheses hold.\\
\textbf{$(H_p)$}
$2<p< \frac{2n-2}{n-2}:=\frac{2^*}{2}+1,n\geq3$, where $2^*$ is the critial cone Sobolev exponent.\\
\textbf{$(H_g)$}
$g(x)\in L^{\infty}\left(int\mathbb{B}\right)\cap C\left(int\mathbb{B}\right) $ is a nonnegative weighted function and
$g(x)\leq \beta$, $x\in int\mathbb{B}$ for positive some constant $\beta>0$.

The main purpose of this paper is to give a systematic investigation on the dynamic behavior of
the solutions to the problem \eqref{1.1} under different initial energy levels.  Before further exploring the problem \eqref{1.1}, we briefly
introduce some topics of the wave type equations associated with this paper.

Research on damping and source terms in wave equations is of great significance for engineering and scientific applications in various fields, as it can enhance system performance and stability, and derive the development and innovation of related technologies. It is well known  that the research of the classical semilinear wave equation on the
bounded domain $\Omega\subset \mathbb{R}^n$
\begin{equation}\label{1.2}
u_{tt}-\Delta u+g(u_t)=f(u)
\end{equation}
has been  extensively studied by  many different unlinked tools. It is impossible to list these results systematically. We only refer to \cite{arrieta1992damped,ball1977remark,levine1997global, liu2020era, pucci1998global,sun2006} and  the references therein. As a special case of \eqref{1.2}, Georgiev and Todorova \cite{georgiev1994existence} conducted the pioneering study on the
interaction between nonlinear damping and polynomial source terms. They studied  the following damped wave
equation
\begin{equation}\label{1.3}
u_{tt}-\Delta u+a u_t|u_t|^{m-1}=b u|u|^{p-1},
\end{equation}
where  $p,m>1$, established  the existence of a global solution to problem \eqref{1.3} when $p\leq m$, and  the blow-up result by constructing
a suitable Lyapunov functional when $p > m$.  Based on their work, Ikehata \cite{ikehata1996some} studied the local existence of the solution and the global existence of the solution under $E(0)<d$. Further more, Vitillaro \cite{vitillaro1999global} investigated the blow-up phenomenon of solutions under  positive initial energy. Recently,
Xu et al \cite{xu22018global} considered the initial boundary value problem of the following fourth-order wave equation  with nonlinear damping and source terms
\begin{equation*}
u_{tt}+\Delta^{2}u+\mu|u_{t}|^{q-2}u_{t}+au=|u|^{p-2}u,
\end{equation*}
and proved  the global existence, energy decay estimate,  blow-up of solution at both $E(0) < d$ and $E(0) = d$  and  blow-up in finite time of solution at arbitrary initial energy $E(0) > 0$ for the case $q=2$ (by the argument as in \cite{gasq2006} ).

Now, we concentrate on  some partial differential equations with potential term. For bounded potentials or potentials with moderate singularities, the unique existence and behavior of solutions to the linear parabolic equation
$u_t-\Delta u-V(|x|)u=0$  is similar with the corresponding properties of the heat equation since the potential $V$ is not big enough. But, such situation may change dramatically for very singular potentials, such as the so-called inverse-square potential $V(x)\sim \lambda |x|^{-2}$ from many physical models (see \cite{singular1995,singular2008}). These singular potentials can be found  in the nonlinear heat \cite{singular1984, singular2008e} and wave equations \cite{singularbook}, nonlinear Schr\"{o}dinger equation  \cite{singular2007f} and the references therein for a detailed discussion of the backgrounds.

In recent years, the analysis on manifolds with conical singularities and the properties of partial differential equations  are intensively studied.  For a singular manifold domain with conical points on the boundary,  it is no longer possible to define the derivative as that in  the classical domain $\Omega \subset \mathbb{R}^n$.  However, the  pseudo-differential operator can reflect singular structure of manifold, such as Fuchsian type Laplace
operator \cite{amann2013function,kondrat1967boundary}, as well as  Schulze and  his collaborators (see \cite{schulze1998boundary,singular1997} and the references therein).

Recently, Chen et al.
\cite{chen2012cone,chen2011existence} have established the corresponding Sobolev inequality and
Poincar\'{e} inequality on the cone Sobolev spaces. From then on, based on those precursory results, there is a number of works on the study of nonlinear partial differential  equations  associated with the manifold with
conical singularities. For example, we can refer the readers  to \cite{chenliu,di2020}  on the initial boundary value problem of the nonlinear parabolic type equations with conical degeneration
on the cone Sobolev spaces. Some results about  well-posedness of the problems were obtained by potential wells theory  proposed by Payne and Sattinger \cite{payne1975saddle,sattinger1968global}. Here we also mention the work \cite{chen2023high}, where  Chen, R\u{a}dulescu and  Xu  \cite{chen2023high}  investigated the following initial boundary value
problem for semilinear conical degenerate parabolic equations with singular potential
\begin{equation*}
u_{t}-\Delta_{\mathbb{B}}u+\kappa V(x)u=g(x)|u|^{p-1}u.
\end{equation*}
They obtained the global existence and nonexistence of solutions for high initial energy and established the finite time blow-up of solutions including the upper bound of blow-up time for arbitrary positive initial energy $J(u_0) > 0$ by applying the Levine's concavity method and the lower
bound of the blow-up time and blow-up rate for arbitrary initial energy.

 Alimohammady et al. \cite{alimohammady2013existence, alimohammady2017invariance, alimohammady2017global} investigated  the
initial boundary value problem for the nonlinear wave equations with conical degeneration
on the cone Sobolev spaces. They established  some results about the global existence, nonexistence and asymptotic behavior of the solutions by using a family of potential well and concavity methods. Especially, Alimohammady and Kalleji \cite{alimohammady2013existence} explored the following initial boundary value problem for a class of semilinear  evolution equation
\begin{equation}\label{1.4}
\begin{cases}
\partial_t^k u-\Delta_{\mathbb{B}}u+ V(x)u=g(x)|u|^{p-1}u ,~~x\in int\mathbb{B},t>0, k\geq1\\[2mm]
 \partial_t^{k-1}u(x, t)=0,\,\,\, \,\,\,\,\qquad \qquad \qquad \qquad x\in \partial\mathbb{B},t\geq0,\\[2mm]
u(x,0)=u_{0}(x), ~~~~~~~~~~~~~~~~~~~\,\,\,\,\,\,\,\,\,\,\,\,\,\,\,\,\, x\in int\mathbb{B}.[2mm]
\end{cases}
\end{equation}
They established the global existence, the asymptotic behavior  and the finite time blow-up for  sub-critical  and critical  initial energy cases. However, there are  some bugs in  the proofs of the theorems in \cite{alimohammady2013existence}. Later, Luo, Xu and Yang \cite{luo2022global} amended the above mistakes existed in \cite{alimohammady2013existence}.  More precisely,  the results in \cite{alimohammady2013existence} are not valid for the case $k\geq2$ and the  parameter $\gamma$ of the singular potential in the model \eqref{1.4} should be added. That is, Luo, Xu and Yang \cite{luo2022global}  considered the  following semilinear hyperbolic equations with singular potentials on the manifolds with conical singularities
\begin{equation}\label{1.5}
\begin{cases}
u_{tt}-\Delta_{\mathbb{B}}u-\gamma V(x)u=g(x)|u|^{p-1}u ,~~x\in int \mathbb{B}, t>0\\[2mm]
 u(x, 0)=u_0(x),\,\,\, u_t(x, 0)=u_1(x), \,\,\,\,\,\,\,\,\,\, x\in int\mathbb{B},\\[2mm]
u(x,t)=0, ~~~~~~~~~~~~~~~~~~~\,\,\,\,\,\,\,\,\,\,\,\,\,\,\,\,\,\,\qquad x\in\partial\mathbb{B} ,t\geq0.
\end{cases}
\end{equation}
they proved the local existence and uniqueness of the solution by using the contraction mapping principle and introduced a family of potential wells to derive a threshold of the existence of global solutions and blow up in finite time of solution in both cases with sub-critical and critical initial
energy. Moreover, two sets of sufficient conditions for initial data leading to blow up result were established at arbitrarily positive initial energy level when $\gamma$ is in some value range. However, there are still  some flaws
in the proof of the local existence and uniqueness of the solution and the the cone Hardy inequality \cite[Lemma 2.5]{luo2022global}.

For the wave equations with conical degeneration and potential term $V(x)u$, we can refer to the references \cite{alimohammady2017invariance,xu2019potential} concerning on linear damping. More precisely, the following  initial-boundary value
problem
\begin{equation}\label{1.6}
  u_{tt}-\Delta_{\mathbb{B}}u+V(x)u+\gamma u_t=g_t(x)|u|^{p-1}u,\quad x\in int\mathbb{B}, t>0,
\end{equation}
was studied.  However, the singular potential term in \eqref{1.6} is the same as that in \eqref{1.4}( some bugs were pointed out in \cite{luo2022global}) and linear damping in \eqref{1.6} should be small, i.e. $\gamma \in [0, (p-1)\sqrt{1+C_*^2}\lambda_1^{\frac{1}{2}}]$.  For more types of equations on conical singularities, we can
refer the readers to \cite{alimohammady2020blow,luo2024,zhang2020} and the references therein.

However, when the nonlinear damping  term $|u_t|^{m-2}u_t$ is added into problem \eqref{1.5}, what will happen for problem \eqref{1.1}?

  Here, for
this model, the appearance of singular potential term, nonlinear damping term and source term cause some difficulties.
Motivated by the above  mentioned literature, we shall conduct a comprehensive and systematic study of problem \eqref{1.1}.

The main results and novelties of our work are summarized as follows.\\
(1) We consider the semilinear wave equation with nonlinear damping $|u_t|^{m-2}u_t$
and singular potential term $\gamma V(x)u$. Different from the classical method, i.e.  the contraction mapping principle (see \cite{georgiev1994existence,liu2020era,xu2018global}), we prove the unique and local  existence of solutions by maximal monotone operator theory modifying the argument of \cite{qin2024quasi,zheng2004nonlinear}.\\
(2)We prove the global existence and the energy decay rate at subcritical and critical initial energy
levels by the potential well theory and Nakao's inequality.  The global existence of the solution is also obtained provided   the damping term  dominates the sourc term, i.e. $m\geq p$. Here, we need some techniques to deal with the singular potential term. \\
(3) We obtain he blow-up in finite time of solution at  subcritical and critical initial energy levels using some ordinary differential inequality.
Moreover, we also prove the blow-up in finite time of solution at the arbitrarily high initial energy and the estimate
of lifespan. As far as we know, there are no results about the finite time blow-up result of solutions to wave-type equation  with
nonlinear damping and singular potential term at arbitrarily high initial energy. The main idea to prove the finite time blow-up of solution of \eqref{1.1} at high energy level
is to modify the method used in recent literature \cite{liao2023blow} in a simply way. We point out that this method would be also effective for the some
other model with nonlinear damping.\\
(4) We also give a necessary and sufficient condition for solutions blowing up in finite time when $2 < p < 2 + \frac{4}{n}$.
 (see Theorem \ref{th5.5}).

The rest of this paper is organized as follows. For the convenience of readers, in Section 2, we first introduce some definitions and properties
of cone Sobolev spaces and the corresponding properties. We also give some lemmas which will  be used  in establishing the main theorems. In Section 3, we are devoted to
prove the local well-posedness of solutions. In Section 4, the global existence and the energy decay
estimate are derived. Finally, the blow-up phenomena of solutions  with   different initial energy levels
are given in Section 5.


\section{Preliminaries}
\setcounter{equation}{0}
\subsection{ Cone Sobolev spaces}
 In this subsection,  we introduce some definitions and properties of cone Sobolev spaces briefly, which is enough to make our paper readable (see \cite{chen2012cone,chen2011existence,schulze1998boundary} for details ).

 Let $X$ be a closed, compact, $C^{\infty}$ manifold, and $X^{\Delta}=(\overline{\mathbb{R}}_+\times X)/(\{0\}\times X)$  be a local
model interpreted as a cone with the base $X$.  We denote by $X^\Lambda= \mathbb{R}_+\times X$ as the corresponding
open stretched cone with $X$. An $n$-dimensional manifold $B$ with conical singularities is a
topological space with a finite subset $B_0=\{b_1,b_2\cdots,b_M\}\subset B$ of conical singularities. For
simplicity, we assume that the manifold $B$ has only one conical point on the boundary. Hence, near the conical point, we have a stretched manifold $\mathbb{B}$, associated with B.
   \begin{definition}(\cite{chen2012cone})
    Let $\mathbb{B} = [0, 1) \times  X$ be a stretched manifold of the manifold $B$ with
	conical singularity. Then the cone Sobolev space $\mathcal{H}_{p}^{m,\gamma}(\mathbb{B})$ for $m\in \mathbb{N}, \gamma\in \mathbb{R},p\in (1, \infty)$ is defined as
	$$\mathcal{H}_{p}^{m,\gamma}(\mathbb{B}):=\left\{W_{loc}^{m,p}(int \mathbb{B})|\omega u\in \mathcal{H}_{p}^{m,\gamma}(X^{\wedge})\right\}$$
	for any cut off function $\omega$ supported by a collar neighborhood of $[0, 1)\times \partial \mathbb{B}$.  Moreover,  the subspace $\mathcal{H}_{p,0}^{m,\gamma}(\mathbb{B})$  of $\mathcal{H}_{p}^{m,\gamma}(\mathbb{B})$is defined  as
	 $$\mathcal{H}_{p,0}^{m,\gamma}(\mathbb{B}):=[\omega]\mathcal{H}_{p,0}^{m,\gamma}(X^{\wedge})+[1-\omega]W_0^{m,p}(int \mathbb{B}), $$
	$W_0^{m,p}(int \mathbb{B})$ denotes the closure  of  $C_0^{\infty}\left(int\mathbb{B}\right)$ in  the Sobolev space $W^{m, p}(\tilde{X})$ where $\tilde{X}$ is a closed compact $C^{\infty}$ manifold of dimension $n$ that containing $\mathbb{B}$ as a submanifold with boundary.
	
	The space  $\mathcal{H}_{p}^{m,\gamma}(\mathbb{B})$ is a Banach space for $p\in[1,\infty)$, a Hilbert space for $p=2$. The cone Lebesgue space $L_p^{\gamma}(\mathbb{B}):=\mathcal{H}_p^{0,\gamma}(\mathbb{B})$. Especially, $L_p^{\frac{n}{p}}(\mathbb{B}):=\mathcal{H}_p^{0,\frac{n}{p}}(\mathbb{B})$ endowed with the norm
	 $$\|u\|_{L_p^{\frac{n}{p}}(\mathbb{B})}=\left(\int_{\mathbb{B}}x_1^n|x_1^{-\frac{n}{p}}u(x)|^p\frac{dx_1}{x_1}dx'\right)^{\frac{1}{p}}=\left(\int_{\mathbb{B}}|u(x)|^{p}\frac{dx_1}{x_1}dx'\right)^{\frac{1}{p}}<\infty.$$
 \end{definition}

 \begin{lemma}\label{a}(Cone Poincar\'{e} inequality\cite{chen2012cone})  Let $1<p<+\infty$ and $\gamma\in \mathbb{R}$. If $u(x)\in {\mathcal{H}}_{p,0}^{1,\gamma}\left(\mathbb{B}\right)$,
 then $$\|u(x)\|_{L^{\gamma}_{p}\left(\mathbb{B}\right)}\leq c\|\nabla_{\mathbb{B}}u(x)\|_{L^{\gamma}_{p}\left(\mathbb{B}\right)} ,$$
 where $\nabla_{\mathbb{B}}=(x_1\partial_{x_1},\partial_{x_2},\cdot\ldots,\partial_{x_n})$ is the gradient operator on $\mathbb{B}$ and the constant $c$
depends only on $\mathbb{B}$ and $p$.
 \end{lemma}

 \begin{lemma}\label{lemmaconeemb}(Cone Sobolev embedding \cite{chen2011existence} ) For $1 \leq q< 2^{*} =\frac{2n}{n-2}$,  the embedding $\mathcal{H}_{2,0}^{1,\frac{n}{2}}(\mathbb{B})\hookrightarrow L_q^{\frac{n}{q}}(\mathbb{B}) $ is continuous and compact.
\end{lemma}

\begin{lemma}\label{b}\cite{luo2022global}
For all $u\in{\mathcal{H}}_{2,0}^{1, \frac{n}{2}}\left(\mathbb{B}\right)\setminus\left\{0\right\}$,the following
inequality $$ \| g ( x ) ^ { \frac { 1 } { p  } } u \| _ { L _ { p  } ^ { \frac { n } { p } }\left(\mathbb{B}\right) } \leq C _ { * } \|\nabla _ {\mathbb{B} } u \| _ { L _ { 2 } ^ { \frac { n } { 2 } } \left(\mathbb{B}\right)  }$$
holds for
$$ C_{*}: =\mathop{sup}\limits_{u\in{\mathcal{H}}_{2,0}^{1, \frac{n}{2}}\left(\mathbb{B}\right)\setminus\left\{0\right\}}
\frac{\|g(x)^\frac{1}{p}u\|_{L_{p}^{\frac{n}{p}}\left(\mathbb{B}\right)}}
{\|\nabla_{\mathbb{B}}u\|^{2}_{L_{2}^{\frac{n}{2}}\left(\mathbb{B}\right)}}.$$
\end{lemma}
\begin{lemma}\label{c}\cite{chen2012cone} There exist $ 0 < \lambda _ {1}<\lambda _ {2}\leq \lambda_{3}
\leq \cdots \leq \lambda _ { k } \leq \cdots \rightarrow
+ \infty$, such that for each $k\geq1$, the following Dirichlet problem
\begin{equation*}
   \begin{cases}
  -\Delta_{\mathbb{B}}\omega_{k}=\lambda_{k}\omega_{k},in\,\, int\mathbb{B},\\[2mm]
 \omega_{k}=0,on\,\, \partial\mathbb{B}, \\[2mm]
  \end{cases}
\end{equation*}
admit a non-trivial solution in ${\mathcal{H}}_{2,0}^{1,\frac{n}{2}}\left(\mathbb{B}\right)$. Moreover,
$\left\{\omega_k\right\}_{k\geq1}$ constitute an orthonormal basis of the Hilbert space
${\mathcal{H}}_{2,0}^{1,\frac{n}{2}}\left(\mathbb{B}\right)$, and the inequality
 $\lambda_1\|u\|^2_{L^{\frac{n}{2}}_{2}\left(\mathbb{B}\right)}\leq
 \|\nabla_{\mathbb{B}}u\|^{2}_{L_{2}^{\frac{n}{2}}\left(\mathbb{B}\right)}$ holds.
\end{lemma}

\begin{lemma}\label{d}(Cone Hardy inequality\cite{chen2012existence})
For all  $u\in H_{2,0}^{1,\frac{n}{2} }(\mathbb{B})\backslash \{0\}$, there exists an optimal positive constant $C^*$ such that
	\begin{equation}\label{2.1}
		\|V(x)^{\frac{1}{2}}u\|_{L^{\frac{n}{2}}_2(\mathbb{B})}\leq C^*\|\nabla_{\mathbb{B}} u\|_{L^{\frac{n}{2}}_2(\mathbb{B})},
	\end{equation}
where the following two  kinds of singular potential functions are given
$$V_1(x_1,x')=\left(\frac{n-3}{2}\right)^2\frac{1}{x_1^2+|x'|^2},\quad V_2(x_1,x')= \left(\frac{n-1}{2}\right)^2\frac{x_1^{-2}e^{-\frac{1}{x_1^2}}}{e^{-\frac{1}{x_1^2}}+|x'|^2}.$$
\end{lemma}

\begin{remark} \label{Re2.1}
Here we need to amend the cone Hardy inequality used in \cite{luo2022global}. For  all $u\in{\mathcal{H}}_{2,0}^{1,\frac{n}{2}}\left(\mathbb{B}\right)\setminus\left\{0\right\}$, one has the  following
inequality
\begin{equation*}
 c_1\|\nabla_{\mathbb{B}}u\|_{L^{\frac{n}{2}}_2(\mathbb{B})}^2\leq \|\nabla_{\mathbb{B}}u\|_{L^{\frac{n}{2}}_2(\mathbb{B})}^2-\gamma \|\sqrt{V(x)}u\|_{L^{\frac{n}{2}}_2(\mathbb{B})}^2\leq c_2\|\nabla_{\mathbb{B}}u\|_{L^{\frac{n}{2}}_2(\mathbb{B})}^2,
 \end{equation*}
 where the positive constants $c_1$ and $c_2$ are selected as follows
 \begin{equation*}
   c_1:=\begin{cases}1-\gamma(C^{*})^2,&\gamma>0;\\
   1,&\gamma\leq0,
   \end{cases}  \text{and } \,\,\,  c_2:=\begin{cases} 1,&\gamma>0;\\
   1-\gamma(C^{*})^2,&\gamma\leq0.
   \end{cases}
 \end{equation*}
\end{remark}

Notice that if  $u(x)\in {L^{\frac{n}{p}}_{p}\left(\mathbb{B}\right)}$
 and $v(x)\in {L^{\frac{n}{q}}_{q}\left(\mathbb{B}\right)}$ for $p, q\in(1,\infty)$,
 $\frac{1}{p}+\frac{1}{q}=1$, then we have the following cone H\"{o}lder inequality
 $$(u, v):=\int_{\mathbb{B}}u(x)v(x)\frac{dx_1}{x_1}dx'\leq \|u(x)\|_{L^{\frac{n}{p}}_{p}\left(\mathbb{B}\right)}\|v(x)\|_{L^{\frac{n}{q}}_{q}\left(\mathbb{B}\right)}.$$

In order to state our main results, we give  the definition of a weak solution of problem \eqref{1.1}.
\begin{definition}(Weak solution) Let $T>0$, $u_0\in \mathcal{H}_{2,0}^{1,\frac{n}{2}}
 \left(\mathbb{B}\right) $ and $u_1\in L_2^{\frac{n}{2}}(\mathbb{B})$.
By a weak solution to problem \eqref{1.1}, we mean a  function $u=u(x,t)$ satisfies
$$u \in C \left([0,T];
\mathcal{H}_{2,0}^{1,\frac{n}{2}}
\left(\mathbb{B}\right)\right)\cap C^1\left([0,T];{L^{\frac{n}{2}}_{2}\left(\mathbb{B}\right)}\right)\cap
C^2 \left([0,T];\mathcal{H}_{2}^{-1,-\frac{n}{2}}\left(\mathbb{B}\right)\right)
$$
such that $u(0)=u_0$, $u_t(0)=u_1$, and
\begin{equation}\label{2.2}
\langle u_{tt},\eta \rangle+
(\nabla_{\mathbb{B}}u,\nabla_{\mathbb{B}}\eta)
-(\gamma V(x)u,\eta)+(|u_t|^{m-2}u_t,\eta)=(g(x)|u|^{p-2}u,\eta)
\end{equation}
holds for any test function $\eta \in \mathcal{H}_{2,0}^{1,\frac{n}{2}}
\left(\mathbb{B}\right)\cap{L^{\frac{n}{m}}_{m}\left(\mathbb{B}\right)}$
and for almost $t \in [0,T]$, where  $\langle \cdot,\cdot \rangle$ denotes
 the duality pairing between $\mathcal{H}_{2,0}^{1,\frac{n}{2}}
 \left(\mathbb{B}\right)$ and $\mathcal{H}_{2}^{-1,-\frac{n}{2}}\left(\mathbb{B}\right)$.
\end{definition}

\subsection{Some auxiliary results}
In this subsection, we give some properties of energy functional and some technical lemmas which will be used to prove the  main results.

Now, we introduce the following functionals  on cone Sobolev space     $\mathcal{H}_{2,0}^{1,\frac{n}{2}}(\mathbb{B})$. The potential energy functional
\begin{equation} \label{2.3}
J(u):=\frac{1}{2}\|\nabla_{\mathbb{B}}u\|^{2}_{L_{2}^{\frac{n}{2}}\left(\mathbb{B}\right)}
-\frac{1}{2}\gamma\|V(x)^{\frac{1}{2}}u\|^{2}_{L_{2}^{\frac{n}{2}}\left(\mathbb{B}\right)}
-\frac{1}{p}\|g(x)^\frac{1}{p}u\|_{L_{p}^{\frac{n}{p}}\left(\mathbb{B}\right)}^{p},
\end{equation}
 the Nehari functional

 \begin{equation} \label{2.4}
I(u):=\|\nabla_{\mathbb{B}}u\|^{2}_{L_{2}^{\frac{n}{2}}\left(\mathbb{B}\right)}
-\gamma\|V(x)^{\frac{1}{2}}u\|^{2}_{L_{2}^{\frac{n}{2}}\left(\mathbb{B}\right)}
-\|g(x)^\frac{1}{p}u\|_{L_{p}^{\frac{n}{p}}\left(\mathbb{B}\right)}^{p},
\end{equation}
and the total energy functional

\begin{equation} \label{2.5}
E(t):=\frac{1}{2}\|u_t\|^2_{L^{\frac{n}{2}}_{2}\left(\mathbb{B}\right)}+J(u).
\end{equation}
The corresponding Nehari manifold can be denoted by
$$\mathcal{N}:=\left\{u\in {\mathcal{H}}_{2,0}^{1,\frac{n}{2}}\left(\mathbb{B}\right)
|I(u)=0,\|\nabla_{\mathbb{B}}u\|_{L_{2}^{\frac{n}{2}}\left(\mathbb{B}\right)}\neq0\right\}$$
and the depth of the potential well or the so-called mountain pass level
\begin{equation}\label{defd}
  d:=\inf\limits_{u\in\mathcal{N}}J(u)=\inf \left\{\sup _{\lambda \geq 0} J(\lambda u), u \in \mathcal{H}_{2,0}^{1, \frac{n}{2}}(\mathbb{B}),\left\|\nabla_{\mathbb{B}} u\right\|_{L_2^{\frac{n}{2}}(\mathbb{B})} \neq 0\right\}.
\end{equation}

Now, we define the potential well
$$\mathcal{W}:=\left\{u\in {\mathcal{H}}_{2,0}^{1,\frac{n}{2}}\left(\mathbb{B}\right)
|I(u)>0,J(u)<d\right\}\cup\left\{0\right\}$$
and the outer of the potential well
$$\mathcal{V}:=\left\{u\in {\mathcal{H}}_{2,0}^{1,\frac{n}{2}}\left(\mathbb{B}\right)
|I(u)<0,J(u)<d\right\}.$$

\begin{lemma}\label{e}(Conservation law of the energy)
Assume that $u(x,t)$ is the  solution of the problem \eqref{1.1}, then  one has
\begin{equation}\label{2.7}
E(t)+\int_{0}^{t}\|u_{\tau}\|^{m}_{L^{\frac{n}{m}}_{m}\left(\mathbb{B}\right)}d\tau=E(0).
\end{equation}\end{lemma}
\begin{proof}
By testing \eqref{1.1} with $u_t$ and integrating with respect to $t$ over $[0, t]$, we have
\begin{equation*}
   \begin{split}&\int_0^t\langle u_{tt},u_{t}\rangle d\tau
   +\int_0^t(\nabla_{\mathbb{B}}u,\nabla u_{t})d\tau
   -\gamma\int_0^t(V(x)u,u_{t})d\tau
   +\int_0^t(|u_t|^{m-2}u_t,u_t)d\tau\\[2mm]
   &=g(x)\int_0^t(|u|^{p-2}u,u_t)d\tau,
   \end{split}
\end{equation*}
that is,
\begin{equation*}
   \begin{split}&\frac{1}{2}\|u_t\|^2_{L^{\frac{n}{2}}_{2}\left(\mathbb{B}\right)}   +\frac{1}{2}\|\nabla_{\mathbb{B}}u\|^{2}_{L_{2}^{\frac{n}{2}}\left(\mathbb{B}\right)}
-\frac{1}{2}\gamma\|V(x)^{\frac{1}{2}}u\|^{2}_{L_{2}^{\frac{n}{2}}\left(\mathbb{B}\right)}
-\frac{1}{p}\|g(x)^\frac{1}{p}u\|_{L_{p}^{\frac{n}{p}}\left(\mathbb{B}\right)}^{p}
+\int_{0}^{t}\|u_{\tau}\|^{m}_{L^{\frac{n}{m}}_{m}\left(\mathbb{B}\right)}d\tau\\[2mm]
   &=\frac{1}{2}\|u_t(0)\|^2_{L^{\frac{n}{2}}_{2}\left(\mathbb{B}\right)}   +\frac{1}{2}\|\nabla_{\mathbb{B}}u(0)\|^{2}_{L_{2}^{\frac{n}{2}}\left(\mathbb{B}\right)}
-\frac{1}{2}\gamma\|V(x)^{\frac{1}{2}}u(0)\|^{2}_{L_{2}^{\frac{n}{2}}\left(\mathbb{B}\right)}
-\frac{1}{p}\|g(x)^\frac{1}{p}u(0)\|_{L_{p}^{\frac{n}{p}}\left(\mathbb{B}\right)}^{p},
   \end{split}
\end{equation*}
which implies $E(t)+\int_{0}^{t}\|u_{\tau}\|^{m}_{L^{\frac{n}{m}}_{m}}d\tau=E(0).$
\end{proof}

For the completeness, we give the following Lemmas (Lemma  \ref{i} and  Lemma \ref{k})  by modifying the proofs in \cite{luo2022global}.
\begin{lemma}\cite{luo2022global}\label{i} Suppse that $u\in {\mathcal{H}}_{2,0}^{1,\frac{n}{2}}\left(\mathbb{B}\right)$
and $\|\nabla_{\mathbb{B}}u\|_{L_{2}^{\frac{n}{2}}\left(\mathbb{B}\right)}\neq0$, we obtain\\
\textbf{(i)}
$\lim\limits_{\lambda\to 0}J(\lambda u)=0$
and $\lim\limits_{\lambda\to +\infty}J(\lambda u)=-\infty$.\\
\textbf{(ii)}
in the interval $(0,\infty)$, there exists a unique $\lambda^*=\lambda^*(u)$ such that $\frac{d}{d\lambda}|_{\lambda=\lambda^*}=0$,  $J(u)$ is increasing on $0\leq\lambda\leq\lambda^*$, decreasing on $\lambda^*\leq\lambda<\infty$ and obtains the maximun at $\lambda=\lambda^*$.\\
\textbf{(iii)}
$I(\lambda u)>0$ for $0\leq\lambda\leq\lambda^*$, $I(\lambda u)<0$ for $\lambda^*\leq\lambda<\infty$, and $I(\lambda^* u)=0$.\\
\textbf{(iv)}
the mountain-pass level $d$ is defined by
$$d:=\mathop{inf}\limits_{u\in\mathcal{N}}J(u)=
\frac{p-2}{2p}C_*^{-\frac{2p}{p-2}}\left(1-\gamma(C^\ast)^{2}\right)^{\frac{p}{p-2}}.$$
\end{lemma}

\begin{lemma}\cite{luo2022global}\label{k}
Assume that $u\in\mathcal{H}_{2,0}^{1,\frac{n}{2}}\left(\mathbb{B}\right)$ and $\gamma\in[0,\frac{1}{{C^*}^2})$.\\
\textbf{(i)}  If $0\leq \|\nabla_{\mathbb{B}}u\|_{L_{2}^{\frac{n}{2}}\left(\mathbb{B}\right)}<
\left(\frac{(1-\gamma(C^\ast)^{2})}{{C_*}^{p}}\right)^{\frac{1}{p-2}}$, then
$I(u)>0$. \\
\textbf{(ii)} If $I(u)<0$, then $\|\nabla_{\mathbb{B}}u\|_{L_{2}^{\frac{n}{2}}\left(\mathbb{B}\right)}>
\left(\frac{(1-\gamma(C^\ast)^{2})}{{C_*}^{p}}\right)^{\frac{1}{p-2}}$.\\
\textbf{(iii)} If $I(u)=0$, then $\|\nabla_{\mathbb{B}}u\|_{L_{2}^{\frac{n}{2}}\left(\mathbb{B}\right)}>
\left(\frac{(1-\gamma(C^\ast)^{2})}{{C_*}^{p}}\right)^{\frac{1}{p-2}}$ or $\|\nabla_{\mathbb{B}}u\|_{L_{2}^{\frac{n}{2}}\left(\mathbb{B}\right)}=0$.
\end{lemma}

\begin{lemma}\label{j}(Invariant sets) Suppose $(H_p)$ and $(H_g)$ hold. Let $u(t), t\in [0, T_{max})$  be the unique solution of problem \eqref{1.1} with initial value $u_0\in {\mathcal{H}}_{2,0}^{1,\frac{n}{2}}\left(\mathbb{B}\right)$ and $u_1\in L^{\frac{n}{2}}_2(\mathbb{B})$, where $T_{max}$ is the maximum existence time of solution.\\
\textbf{(i)} If there exists a $t_0\in [0, T_{max})$ such that $E(t_0)<d$,
$u(t_0)\in \mathcal{W}$, then $u(t)\in \mathcal{W}$ for all $t \in [t_0, T_{max})$. \\
\textbf{(ii)} If there exists a $t_0\in [0, T_{max})$ such that either $E(t_0)<d$ or $E(t_0)=d$ with $(u(t_0), u_t(t_0))>0$, then $u(t)\in \mathcal{V}$ for all  $t \in [t_0, T_{max})$ provided that $u(t_0)\in \mathcal{V}$.\\
\textbf{(iii)} If there exists a $t_0\in [0, T_{max})$ such that $E(t_0)<0$ or $E(t_0)=0$ with $u(t_0)\neq0$, then  $u(t)\in \mathcal{V}$ for all $t \in [t_0, T_{max})$.
\end{lemma}
\begin{proof}
\textbf{(i)} Since $E(t)$ is non-increasing with respect to $t$, one can derive from \eqref{2.5} that
\begin{equation}\label{2.8}
J(u(t))=E(t)-\frac{1}{2}\|u_t\|^2_{L^{\frac{n}{2}}_{2}\left(\mathbb{B}\right)}\leq E(t_0)<d,\,\,t_0\leq t<T_{max}.
\end{equation}

Now, we shall split the proof into two cases by the similar argument as \cite{ZhouAMO2023}.

Case 1: \emph{there exists a $t_1\in [t_0, T_{max})$ such that $I(u(t_1))>0$}.  If the conclusion of part (i) is not true, by the continuity of $I(u(t))$ and $J(u(t))$ with respect to $t$,
 there exists a $t_2\in (t_1, T_{max})$ such that $I(u(t))>0$  for all $t\in [t_1, t_2)$, $u(t_2)\neq0$ (otherwise $u(t_2)\in \mathcal{W}$) and $I(u(t_2))=0$, that is $u(t_2)\in \mathcal{N}$.
 Then by the definition of $d$ (see \eqref{defd}), we have $J(u(t_2))\geq d$, which is a contradiction with \eqref{2.8}.

Case 2: \emph{$I(u(t))\leq 0$ for  all $t\in [t_0, T_{max}).$}  Noticing $u(t_0)\in \mathcal{W}$, we must have $u(t_0)=0$, that is $I(u(t_0))=0$.  Now, we end the proof of Case 2 into the following  two subcases.

Subcase 1: \emph{$I(u(t))=0$ for all $t\in (t_0, T_{max})$.} Then, we must have $u(t)=0 $ for all $t\in [t_0, T_{max})$, which yields that $u(t)\in \mathcal{W}$ for all  $t\in [t_0, T_{max})$. Otherwise, if there exists $t_1\in (t_0, T_{max})$ such that $u(t_1)\neq0$, that is $u(t_1)\in\mathcal{N}$. Then by the definition of $d$ (see \eqref{defd}) again, we have $J(u(t_1))\geq d $,   which is a contradiction with \eqref{2.8}.

Subcase 2: \emph{there exists a $t_1\in (t_0, T_{max})$ such that $I(u(t_1))<0$}. We shall prove this subcase can not occur by contradiction. Indeed, if subcase 2 holds, noticing that $I(u(t_0))=0$, we have $t_2:=\{t\in [t_0, t_1), I(u(t))=0\}$ exists and $t_2<t_1$, which yields that $I(u(t_2))=0$ and  there exists a decreasing sequence $\{t_n\}_{n=1}^{\infty}\subset (t_2, t_1)$ such that
\begin{equation*}
  I(u(t_n))<0,\,\,\,\,\lim_{n\to \infty}t_n=t_2.
\end{equation*}
Thus, it follows from Lemma \ref{a} and  Lemma \ref{b} that
\begin{equation*}
   \begin{split}\|\nabla_{\mathbb{B}}u(t_n)\|^{2}_{L_{2}^{\frac{n}{2}}\left(\mathbb{B}\right)}
   &<\|g(x)^\frac{1}{p}u(t_n)\|_{L_{p}^{\frac{n}{p}}\left(\mathbb{B}\right)}^{p}
   +\gamma\|V(x)^{\frac{1}{2}}u(t_n)\|^{2}_{L_{2}^{\frac{n}{2}}\left(\mathbb{B}\right)}\\[2mm]
   &\leq C_\ast^{p}
   \|\nabla_{\mathbb{B}}u(t_n)\|^{p-2}_{L_{2}^{\frac{n}{2}}\left(\mathbb{B}\right)}
   \|\nabla_{\mathbb{B}}u(t_n)\|^{2}_{L_{2}^{\frac{n}{2}}\left(\mathbb{B}\right)}
   +\gamma(C^\ast)^{2}\|\nabla_{\mathbb{B}}u(t_n)\|^{2}_{L_{2}^{\frac{n}{2}}\left(\mathbb{B}\right)},
   \end{split}
\end{equation*}
which implies $$\|\nabla_{\mathbb{B}}u(t_n)\|_{L_{2}^{\frac{n}{2}}\left(\mathbb{B}\right)}>
\left(\frac{(1-\gamma(C^\ast)^{2})}{{C_*}^{p}}\right)^{\frac{1}{p-2}}.$$
Since $u\in C([0, T_{max}); {\mathcal{H}}_{2,0}^{1,\frac{n}{2}}\left(\mathbb{B}\right))$ and $t_n\rightarrow t_2$ as $n\to \infty, $
letting $n\to \infty$ in the above inequality, we have
$$\|\nabla_{\mathbb{B}}u(t_2)\|_{L_{2}^{\frac{n}{2}}\left(\mathbb{B}\right)}\geq
\left(\frac{(1-\gamma(C^\ast)^{2})}{{C_*}^{p}}\right)^{\frac{1}{p-2}}.$$
Hence, the above analysis yields that $I(u(t_2))=0$ and $u(t_2)\neq0$, i.e. $u(t_2)\in \mathcal{N}$. Then by the definition of $d$ (see \eqref{defd}) again, we have $J(u(t_2))\geq d $,   which is a contradiction with \eqref{2.8}.

\textbf{(ii)} First, we prove  the case $E(t_0)<d$ by contraction. Indeed, if the conclusion is not true, since $u\in C([0, T_{max}); {\mathcal{H}}_{2,0}^{1,\frac{n}{2}}\left(\mathbb{B}\right))$, we assume that there exists a first time $t_1\in [t_0, T_{max})$ such that $u(t_1)\in \partial \mathcal{V}$, where $\partial \mathcal{V}$ denotes the boundary of domain $\mathcal{V}$. From the definition of $\mathcal{V}$ and the
continuity of $J(u)$ and $I(u)$ with respect to $t$, we have either
\begin{equation*}
  I(u(t_1))=0,\,\,\|\nabla_{\mathbb{B}}u\|_{L_{2}^{\frac{n}{2}}\left(\mathbb{B}\right)}\neq0\,\,or \,\,\,J(u(t_1))=d.
\end{equation*}
From \eqref{2.8}, we can see that $J(u(t_1))=d$ is impossible.  And if $I(u(t_1))=0, \|\nabla_{\mathbb{B}}u\|_{L_{2}^{\frac{n}{2}}\left(\mathbb{B}\right)}\neq0$, then from \eqref{defd} it is clear that $J(u(t_1))\geq d$, which contradicts \eqref{2.8} again.

Second, we prove the case  $E(t_0)=d$ with $(u(t_0), u_t(t_0))>0$ by contraction. Indeed, if the conclusion is not true, it follows from $u(t_0)\in \mathcal{V}$ and $u\in C([0, T_{max}); {\mathcal{H}}_{2,0}^{1,\frac{n}{2}}\left(\mathbb{B}\right))$ that there must exist a $t_1\in (t_0, T_{max})$ such that $u(t)\in \mathcal{V}$ for all $t\in [t_0, t_1)$ and (a) $J(u(t_1))<d, I(u(t_1))=0$ or (b) $J(u(t_1))=d, I(u(t_1))\leq 0$.  Since $I(u(t))<0$ for all $t\in [t_0, t_1)$, then it follows from Lemma \ref{k} that $$\|\nabla_{\mathbb{B}}u\|_{L_{2}^{\frac{n}{2}}\left(\mathbb{B}\right)}>
\left(\frac{(1-\gamma(C^\ast)^{2})}{{C_*}^{p}}\right)^{\frac{1}{p-2}},\,\,\,t\in [t_0, t_1)$$ which implies $$\|\nabla_{\mathbb{B}}u(t_1)\|_{L_{2}^{\frac{n}{2}}\left(\mathbb{B}\right)}\geq
\left(\frac{(1-\gamma(C^\ast)^{2})}{{C_*}^{p}}\right)^{\frac{1}{p-2}}>0.$$
  From the above inequality, if (a) occurs, we get $u(t_1)\in \mathcal{N}$.  Then by the definition of $d$ (see \eqref{defd}) again, we have $J(u(t_1))\geq d $,   which is a contradiction with \eqref{2.8}.
  If (b) occurs, by Lemma \ref{e},  that is
    \begin{equation*}
    \|u_t(t_1)\|^2_{L_{2}^{\frac{n}{2}}\left(\mathbb{B}\right)}+J(u(t_1))+\int_{t_0}^{t_1}\|u_t(t)\|^m_{L_{m}^{\frac{n}{m}}\left(\mathbb{B}\right)}dt=E(t_0)=d.
  \end{equation*}
Then combining with $J(u(t_1))=d$, we have
 \begin{equation}\label{2.9}
    \|u_t(t_1)\|^2_{L_{2}^{\frac{n}{2}}\left(\mathbb{B}\right)}+\int_{t_0}^{t_1}\|u_t(t)\|^m_{L_{m}^{\frac{n}{m}}\left(\mathbb{B}\right)}dt=0.
  \end{equation}
  On other hand, it is clear from $(u(t_0), u_t(t_0))>0$ that $\|u_t(t_0)\|_{L_{m}^{\frac{n}{m}}\left(\mathbb{B}\right)}\neq0$  by cone H\"{o}lder inequality, which implies that $\int_{t_0}^{t_1}\|u_t(t)\|^m_{L_{m}^{\frac{n}{m}}\left(\mathbb{B}\right)}dt>0$.  Hence, the equality \eqref{2.9} does not hold.   Hence we complete the proof of (ii).

\textbf{(iii)} Noticing \eqref{2.3}-\eqref{2.5}, we have
\begin{equation*}
  E(t)=\frac{1}{2}\|u_t(t)\|^2_{L_{2}^{\frac{n}{2}}\left(\mathbb{B}\right)}+\frac{1}{p}I(u)
  +\frac{p-2}{2p}
  (\|\nabla_{\mathbb{B}}u\|^{2}_{L_{2}^{\frac{n}{2}}\left(\mathbb{B}\right)}
  -\gamma\|V(x)^{\frac{1}{2}}u\|^{2}_{L_{2}^{\frac{n}{2}}\left(\mathbb{B}\right)}).
\end{equation*}
If $E(t_0)<0$ or $E(t_0)=0$ with $u(t_0)\neq0$, we obviously have $E(t_0)<d$ and
$J(u(t_0)\leq E(t_0)<d$, $I(u(t_0))\leq pE(t_0)-\frac{p-2}{2}
  (\|\nabla_{\mathbb{B}}u\|^{2}_{L_{2}^{\frac{n}{2}}\left(\mathbb{B}\right)}
  -\gamma\|V(x)^{\frac{1}{2}}u\|^{2}_{L_{2}^{\frac{n}{2}}\left(\mathbb{B}\right)})<0$, which implies that $u(t_0)\in \mathcal{V}$. Then the conclusion    is true from part (ii).
\end{proof}

\begin{lemma}\label{m}(Nakao's inequality \cite{nakao1978difference})
Let $\phi(t)$ be a nonincreasing nonnegative function on $[0,\infty)$ and satisfy
$$\phi^{1+r}(t)\leq k_0\left(\phi(t)-\phi(t+1)\right),t\in[0,T]$$
where $k_0$ is a postive constant, $r$ is a nonnegative constant. that

\textbf{(i)}
if $r>0$, then $\phi(t)\leq\left(\phi^{-r}(t)+k_0 r[t-1]^{+}\right)^{-\frac{1}{r}}$,

\textbf{(ii)}
if $r=0$, then $\phi(t)\leq\phi(0)e^{-k_{1}[t-1]^{+}}$, where
$$[t-1]^{+}=max\left\{t-1,0\right\},k_1=log\left(\frac{k_0}{k_0-1}\right).$$
\end{lemma}

\begin{lemma}\label{GN}(Cone Gagliardo-Nirenberg inequality)
Let $1\leq s_1<s_2\leq 2^*=\frac{2n}{n-2}$,
Then for any $u\in \mathcal{H}_{2,0}^{1,\frac{n}{2}}\left(\mathbb{B}\right)$, there exists a positive constant $C>0$ such that
\begin{equation}\label{2.10}
    \|u\|_{L^{\frac{n}{s_2}}_{s_2}\left(\mathbb{B}\right)}\leq
     C \|\nabla_{\mathbb{B}}u\|^{\theta}_{L_{2}^{\frac{n}{2}}\left(\mathbb{B}\right)}
     \|u\|^{1-\theta}_{L_{s_1}^{\frac{n}{s_1}}\left(\mathbb{B}\right)},
\end{equation}
where $\theta\in (0, 1]$ satisfies $\frac1{s_2}=\frac{1-\theta}{s_1}+\theta\left(\frac1 2-\frac{1}{n}\right)$.
\end{lemma}

\begin{proof}
Using H\"{o}lder inequality, we have that
\begin{equation}\label{2.11}
    \begin{aligned}
\int_{\mathbb{B}}|u|^{s_2}\frac{\mathrm{d}x_{1}}{x_{1}}\mathrm{d}x^{\prime}& =\int_\mathbb{B}|u|^{s_2\theta} |u|^{s_2(1-\theta)}\frac{\mathrm{d}x_1}{x_1}\mathrm{d}x^{\prime} \\
&\leq\left(\int_{\mathbb{B}}|u|^{2^*}\frac{\mathrm{d}x_1}{x_1}\mathrm{d}x^{\prime}\right)^{\frac{s_2\theta}{2^*}}
\left(\int_{\mathbb{B}}|u|^{s_1}\frac{\mathrm{d}x_1}{x_1}\mathrm{d}x^{\prime}\right)^{\frac{s_2(1-\theta)}{s_1}},
\end{aligned}
\end{equation}
where $\theta$ satisfies $\frac{s_2(1-\theta)}{s_1}+\frac{s_2\theta}{2^{*}}=1$, that is $\frac1{s_2}=\frac{1-\theta}{s_1}+\theta\left(\frac{1} {2}-\frac{1}{n}\right)$.
Then by cone  Sobolev embedding inequality (see Lemma \ref{lemmaconeemb} )
$$\|u\|_{L^{\frac{n}{2^*}}_{2^*}(\mathbb{B})}\leq \mathcal{C}_*\|\nabla_{\mathbb{B}}u\|_{L^{\frac{n}{2}}_{2}(\mathbb{B})},$$
where $\mathcal{C}_*$ is the best constant of the Sobolev embedding, we can obtain \eqref{2.10} directly from \eqref{2.11}.

\end{proof}

\section{Local existence}
\setcounter{equation}{0}
In this section, we shall give the local well-posedness for the solution of the problem \eqref{1.1} by nonlinear semigroup theorem.
Fristly, we consider the following Hilbert space
\begin{equation} \label{3.1}
\mathcal{H}=\mathcal{H}_{2,0}^{1,\frac{n}{2}}\left(\mathbb{B}\right)
\times L^{\frac{n}{2}}_{2}\left(\mathbb{B}\right)
\end{equation}
equipped with the inner product and normal
\begin{equation} \label{3.2}
(U,\tilde{U})_\mathcal{H}=(u_t,\tilde{u}_t)
+(\nabla_{\mathbb{B}}u,\nabla_\mathbb{B}\tilde{u})
-\gamma(V(x)^{\frac{1}{2}}u,V(x)^{\frac{1}{2}}\tilde{u}),
\end{equation}
and
\begin{equation} \label{3.3}
\|U\|^2_\mathcal{H}=\|u_t\|^2+\|\nabla_{\mathbb{B}}u\|^2-\gamma\|V(x)^{\frac{1}{2}}u\|^2,
\end{equation}
for any $U=(u,u_t)$ and $\tilde{U}=(\tilde{u},\tilde{u}_t)$ in $\mathcal{H}$. Then, the problem \eqref{1.1} can be written as an equivalent Cauchy problem ODE
\begin{equation} \label{3.4}
\begin{cases}
\frac{d U}{dt}+\mathbb{A}U=\mathbb{F}(U),\\[2mm]
U(0)=U_0=(u_0,u_1)\in \mathcal{H},
\end{cases}
\end{equation}
where $$U=(u,u_t)\in \mathcal{H}$$
and $\mathbb{A}:\mathcal{D}(\mathbb{A})\subset \mathcal{H} \to \mathcal{H}$ is the nonlinear operator defined by
$$\mathbb{A}U=\left(\begin{array}{c}-u_t\\-\Delta_{\mathbb{B}} u-\gamma V(x) u+|u_t|^{m-2}u_t\end{array}\right)$$
with the domain
$$\mathcal{D}(\mathbb{A})=\left\{U=(u,u_t)\in \mathcal{H}:u_t\in \mathcal{H}_{2,0}^{1,\frac{n}{2}}\left(\mathbb{B}\right),-\Delta_{\mathbb{B}} u-\gamma V(x) u+|u_t|^{m-2}u_t\in   L^{\frac{n}{2}}_{2}\left(\mathbb{B}\right)\right\}.$$
The nonlinear function $\mathbb{F}$ is defined by
$$\mathbb{F}(U)=\left(\begin{array}{c}0\\g(x)|u|^{p-2}u\end{array}\right).$$
\begin{theorem}\label{zz}
Let $U_0 \in \mathcal{H}$, then there exists a positive constant $T$ depending only on $\|U_0\|_{\mathcal{H}}$ such that problem \eqref{3.4} admits a unique mild solution $U(t)\in C([0,T];\mathcal{H})$ with $U(0)=U_{0}$, i.e., $U(t)$ satisfies the following integral equation
\begin{equation}
 U(t)=e^{\mathbb{A}t}U_0+\int_0^t e^{\mathbb{A}(t-s)}\mathbb{F}(U(s))ds,t\in[0,T].
 \end{equation}
Furthermore, the solution $U$ can be extended to a maximal solution in $[0,T_{max})$ such that either \\
(i) $T_{max}=\infty$, i.e., the problem \eqref{3.4} admits a global mild solution; or\\
(ii) $T_{max}<\infty$, and $$\lim_{t\to T_{max}^{-}}\|U(t)\|_{\mathcal{H}}=\infty, $$
i.e., the mild solution blows up at a finite time $T_{max}$.

Moreover, if $U_0\in D(\mathbb{A}):={\mathcal{H}}_{1}$, then $U(t)\in C([0,T_{max}); {\mathcal{H}}_{1})\cap C^1([0,T_{max}); \mathcal{H})$
is the classical solution of problem \eqref{3.4}.  If $U^1(t)$ and $U^2(t)$ are two mild solutions of problem \eqref{3.4}, then there exists a
positive constant $C_0 = C(U^1(0), U^2(0))$, such that
$$\|U^1(t)-U^2(t)\|_\mathcal{H}\leq e^{C_0T}\|U^1(0)-U^2(0)\|_\mathcal{H},$$
where $C_0=C(U^1(0),U^2(0))$.
\end{theorem}

\begin{proof}
We divide this proof into two steps:

\textbf{Step 1: The operator $\mathbb{A}$ generates a $C_0$-semigroup of contraction $S(t) = e^{\mathbb{A}t}$ on $\mathcal{H}$.}\\
We first prove that $\mathbb{A}$ is maximal monotone. To this aim, we denote $\mathbb{A}$ by the
following sum of two operators
$$\mathbb{A}=\mathbb{P}+\mathbb{Q}$$
where $\mathbb{P}:\mathcal{D}(\mathbb{P})\subset \mathcal{H} \to \mathcal{H}$ and
$\mathbb{P}:\mathcal{D}(\mathbb{P})\subset \mathcal{H} \to \mathcal{H}$ are defined by for any $U\in \mathcal{D}(\mathbb{A})$
$$\mathbb{P}U:=\begin{pmatrix}-u_t\\-\Delta u-\gamma V(x)u\end{pmatrix},\quad\mathbb{Q}U:=\begin{pmatrix}0\\|u_t|^{m-2}u_t\end{pmatrix}$$
with domain
$$\mathcal{D}(\mathbb{P}):=\left(\mathcal{H}_2^{2, \frac{n}{2}}\left(\mathbb{B}\right)\cap \mathcal{H}_{2,0}^{1,\frac{n}{2}}\left(\mathbb{B}\right)\right)\times \mathcal{H}_{2,0}^{1,\frac{n}{2}}\left(\mathbb{B}\right)$$
$$\mathcal{D}(\mathbb{Q}):=\mathcal{H}_{2,0}^{1,\frac{n}{2}}\left(\mathbb{B}\right)\times
L^{\frac{n}{m}}_{m}\left(\mathbb{B}\right)$$
Now, we shall prove that $\mathbb{P}$ is  a maximal monotone operator. Indeed, by equation \eqref{1.1}, it follows from
integration by parts and boundary conditions that for all $U,\tilde{U}\in \mathcal{D}(\mathbb{P}) $
\begin{equation} \label{3.5}
   \begin{split}
   &(\mathbb{P}U-\mathbb{P}\tilde{U},U-\tilde{U})\\[2mm]
   &=(\nabla_{\mathbb{B}}(\tilde{u_t}-u_t),\nabla_{\mathbb{B}}(u-\tilde{u}))
   +(\nabla_{\mathbb{B}}(u-\tilde{u}),\nabla_{\mathbb{B}}(u_t-\tilde{u}_t)\\[2mm]
   &-\gamma(V(x)(\tilde{u}_t-u_t),u-\tilde{u})
   -\gamma(V(x)(u-\tilde{u}),u_t-\tilde{u}_t)\\[2mm]
   &=0,
   \end{split}
\end{equation}
where $U=(u,u_t)$, $\tilde{U}=(\tilde{u},\tilde{u}_t)$, which implies that  $\mathbb{P}$ is monotone.
Then, we prove that $\mathbb{P}$ is maximal monotone, it suffices to prove $ R(I+\mathbb{P})=\mathcal{H}$. To this end, for any given $U^*=(u^*,u_t^*)\in\mathcal{H}$, we must show that
\begin{equation}\label{3.6}
U+\mathbb{P}U=U^*,
\end{equation}
admits  a solution $U=(u,u_t)\in \mathcal{D}(\mathbb{P})$. Equation \eqref{3.6} can be rewritten as the following system
\begin{equation}\label{3.7}
\begin{cases}u-u_t=u^*,\\
u_t-\Delta_{\mathbb{B}} u-\gamma V(x)u=u_t^*.
\end{cases}
\end{equation}
Obviously, \eqref{3.7} can be transferred to the following equation
\begin{equation}\label{3.8}
u-\Delta_{\mathbb{B}} u-\gamma V(x)u=u^*+u_t^*.
\end{equation}
Then  the equation \eqref{3.8} is equivalent to the variational problem
$$\mathcal{B}(u, \tilde{u} )=\mathcal{L}(\tilde{u}),$$
where the bilinear form $\mathcal{B}:\mathcal{H}_{2,0}^{1,\frac{n}{2}}\left(\mathbb{B}\right)\times \mathcal{H}_{2,0}^{1,\frac{n}{2}}\left(\mathbb{B}\right)\rightarrow \mathbb{R}$ is given by
$$\mathcal{B}(u,\tilde{u})=\int_{\mathbb{B}}(u\tilde{u}+\nabla_{\mathbb{B}}u\nabla_{\mathbb{B}}\tilde{u}
-\gamma V(x)u\tilde{u})\frac{dx_{1}}{x_{1}}dx',$$
and  the linear form $\mathcal{L}: \mathcal{H}_{2,0}^{1,\frac{n}{2}}\left(\mathbb{B}\right)\rightarrow \mathbb{R}$ is given by
$$\mathcal{L}(\tilde{u})=\int_{\mathbb{B}}\left(u_t^*+u^*\right)\tilde{u}\frac{dx_{1}}{x_{1}}dx'.$$

We can derive  that  $\mathcal{B}: \mathcal{H}_{2,0}^{1,\frac{n}{2}}\left(\mathbb{B}\right)\times \mathcal{H}_{2,0}^{1,\frac{n}{2}}\left(\mathbb{B}\right)\rightarrow \mathbb{R}$ is a bilinear and continuous form and $\mathcal{L}: \mathcal{H}_{2,0}^{1,\frac{n}{2}}\left(\mathbb{B}\right)\rightarrow \mathbb{R}$ is a linear and continuous form. In fact, for any $u, \tilde{u}\in \mathcal{H}_{2,0}^{1,\frac{n}{2}}\left(\mathbb{B}\right)$,  noticing Lemma \ref{c}, we have
\begin{equation*}
   \begin{split}
|\mathcal{B}(u, \tilde{u})|&=\int_\mathbb{B}u\tilde{u}\frac{dx_{1}}{x_{1}}dx'
+\int_\mathbb{B}\nabla_{\mathbb{B}}u\nabla_{\mathbb{B}}\tilde{u}\frac{dx_{1}}{x_{1}}dx'
-\gamma\int_\mathbb{B}V(x)u\tilde{u}\frac{dx_{1}}{x_{1}}dx'\\[2mm]
&\leq \|u\|_{L^{\frac{n}{2}}_{2}\left(\mathbb{B}\right)}
\|\tilde{u}\|_{L^{\frac{n}{2}}_{2}\left(\mathbb{B}\right)}
+\|\nabla_{\mathbb{B}}u\|_{L^{\frac{n}{2}}_{2}\left(\mathbb{B}\right)}
\|\nabla_{\mathbb{B}}\tilde{u}\|_{L^{\frac{n}{2}}_{2}\left(\mathbb{B}\right)}
-\gamma\|V^{\frac{1}{2}}(x)u\|_{L^{\frac{n}{2}}_{2}\left(\mathbb{B}\right)}
\|V^{\frac{1}{2}}(x)\tilde{u}\|_{L^{\frac{n}{2}}_{2}\left(\mathbb{B}\right)}\\[2mm]
&\leq C\|\nabla_{\mathbb{B}} u\|_{L^{\frac{n}{2}}_{2}\left(\mathbb{B}\right)}\|\nabla_{\mathbb{B}}\tilde{u}\|_{L^{\frac{n}{2}}_{2}\left(\mathbb{B}\right)},
\end{split}
\end{equation*}
and
\begin{equation*}
   \begin{split}
\mathcal{L}(\tilde{u})&=\int_\mathbb{B}(u_t^*+u^*)\tilde{u}\frac{dx_{1}}{x_{1}}dx'\\[2mm]
&\leq \|u_t^*\|_{L^{\frac{n}{2}}_{2}\left(\mathbb{B}\right)}
\|\tilde{u}\|_{L^{\frac{n}{2}}_{2}\left(\mathbb{B}\right)}
+\|u^*\|_{L^{\frac{n}{2}}_{2}\left(\mathbb{B}\right)}
\|\tilde{u}\|_{L^{\frac{n}{2}}_{2}\left(\mathbb{B}\right)}\\[2mm]
&\leq C\|\nabla_{\mathbb{B}}\tilde{u}\|_{L^{\frac{n}{2}}_{2}\left(\mathbb{B}\right)}.
\end{split}
\end{equation*}

Moreover, $\mathcal{B}$ is also  coercive since
\begin{equation*}
   \begin{split}
\mathcal{B}(u, u)= \|u\|^2_{L^{\frac{n}{2}}_{2}\left(\mathbb{B}\right)}
+\|\nabla_{\mathbb{B}}u\|^2_{L^{\frac{n}{2}}_{2}\left(\mathbb{B}\right)}
-\gamma\|V^{\frac{1}{2}}(x)u\|^2_{L^{\frac{n}{2}}_{2}\left(\mathbb{B}\right)}\geq \left(1-\gamma (C^*)^2\right)\|\nabla_{\mathbb{B}} u\|_{L^{\frac{n}{2}}_{2}\left(\mathbb{B}\right)},
\end{split}
\end{equation*}
where we have used cone Hardy inequality Lemma \ref{d}.

Hence, using the Lax-Milgram Theorem \cite{zheng2004nonlinear}, the elliptic equation \eqref{3.8} possesses a unique
weak solution $u\in \mathcal{H}_{2,0}^{1,\frac{n}{2}}\left(\mathbb{B}\right) $. In addition, it follows from system \eqref{3.7} that
\begin{equation*}
\begin{gathered}
u_{t} =u-u^*\in \mathcal{H}_{2,0}^{1,\frac{n}{2}}\left(\mathbb{B}\right), \\
\Delta u =u_t-\gamma V(x)u-u_t^*\in L^{\frac{n}{2}}_{2}\left(\mathbb{B}\right).
\end{gathered}
\end{equation*}
Therefore, we we have shown that there exists $U=(u,u_t)\in\mathcal{D}(\mathbb{P})$ such that  $(I+\mathbb{P})U=U^*$. Since $U^*$ is arbitrary, $R(I+\mathbb{P}) = \mathcal{H}$ is
established. The proof of the maximal
monotonicity of $\mathbb{P}$ is complete.

Next, we shall show that $\mathbb{Q}$ is also a maximal monotone operator. By Hille-Yosida's theorem \cite[Chap.2.2]{zheng2004nonlinear}, it suffices to show that $\mathbb{Q}$ is monotone and hemicontinuous. Indeed, for any $U=(u,u_t),\tilde{U}=(u,\tilde{u})\in \mathcal{H}$, we
have
\begin{equation*}
 \begin{split}
&\left\langle\mathbb{Q}U-\mathbb{Q}\tilde{U},U-\tilde{U}\right\rangle\\[2mm]
&=\int_\mathbb{B}\left(|u_t|^{m-2}u_t-|\tilde{u}_t|^{m-2}\tilde{u}_t\right)(u_t-\tilde{u}_t)\frac{dx_{1}}{x_{1}}dx'\\[2mm]
&\geq C\|u_t-\tilde{u}_t\|^m_{L^{\frac{n}{m}}_{m}\left(\mathbb{B}\right)}\geq 0,
 \end{split}
\end{equation*}
which yields that $\mathbb{Q}$ is monotone in $\mathcal{H}$. Then we will prove
that $\mathbb{Q}$ is hemicontinuous. In fact, for any given $U=(u,u_t),U_i=(u_i,u_{it})\in \mathcal{H}(i=1,2)$, we obviously have
\begin{equation*}
 \begin{split}&
\left|\left\langle\mathbb{Q}(U_1+\lambda U_2),U\right\rangle-\left\langle\mathbb{Q}U_1,U\right\rangle\right|\\[2mm]
=&\int_{\mathbb{B}}\left(|u_{1t}+\lambda u_{2t}|^{m-2}(u_{1t}+\lambda u_{2t})-|u_{1t}|^{m-2}u_{1t}, u_t\right)\frac{dx_{1}}{x_{1}}dx'\\[2mm]
\to &\,\, 0,\,\,\text{as}\,\,\, \lambda \to 0,\\[2mm]
 \end{split}
\end{equation*}
which implies that $\mathbb{Q}$ is  hemicontinuous and thus $\mathbb{Q}$ is maximal monotone.

It thus follows that $\mathbb{P}$ and $\mathbb{Q}$ are both maximal monotone. Furthermore, by the
corollary of Hille-Yosida's Theorem \cite[Chapter 2.2]{zheng2004nonlinear} and int$\mathcal{D}(\mathbb{P})\cap\mathcal{D}(\mathbb{Q})\neq \varnothing$, we conclude that $\mathbb{A} = \mathbb{P} + \mathbb{Q}$ is maximal monotone. Moreover, $\mathcal{D}(\mathbb{A})$ is densely defined in $\mathcal{H}$ and thus the operator $\mathbb{A}$ generates a $C_0$-semigroup by Lumer-Phillips Theorem \cite{barbu2010nonlinear}.

\textbf{Step 2: $\mathbb{F}(U)$ is locally Lipschitz in $\mathcal{H}$.}\\
To this end,  let $\mathcal{B}$ be a bounded set of $\mathcal{H}$ and
$$U_1=(u_1(t),u_{1t}(t)),U_2=(u_2(t),u_{2t}(t))\in\mathcal{B}.$$  From the cone Sobolev embedding,  H\"{o}lder inequality and \eqref{3.3}, we deduce that
\begin{equation}\label{3.10}
   \begin{split}
   &\|f(u_1)-f(u_2)\|_{L^{\frac{n}{2}}_{2}\left(\mathbb{B}\right)}\\[2mm]
   &=\left\{\int_{\mathbb{B}}
   |g(x)|u_1|^{p-2}u_1-g(x)|u_2|^{p-2}u_2|^2\frac{dx_{1}}{x_{1}}dx'\right\}^{\frac{1}{2}}\\[2mm]
   &\leq\beta\left\{\int_{\mathbb{B}}||u_1|^{p-2}u_1-|u_2|^{p-2}u_2|^2\frac{dx_{1}}{x_{1}}dx'\right\}^{\frac{1}{2}}\\[2mm]
  &\leq\beta(p-1)
  \left\{\int_{\mathbb{B}}(|u_1|+|u_2|)^{2(p-2)}|u_1-u_2|^{2}\frac{dx_{1}}{x_{1}}dx'\right\}^{\frac{1}{2}}\\[2mm]  &\leq\beta(p-1)
  \left\{\int_{\mathbb{B}}||u_1|+|u_2||^{(p-2)n}\frac{dx_{1}}{x_{1}}dx'\right\}^{\frac{1}{n}}
  \left\{\int_{\mathbb{B}}|u_1-u_2|^{\frac{2n}{n-2}}\frac{dx_{1}}{x_{1}}dx'\right\}^{\frac{n-2}{2n}}\\[2mm]
&\leq\beta(p-1)
  \|(|u_1|+|u_2|)\|^{p-2}_{L^{\frac{1}{p-2}}_{(p-2)n}\left(\mathbb{B}\right)}
   \|{|\nabla_{\mathbb{B}}(u_1-u_2)}\|_{L^{\frac{n}{2}}_{2}\left(\mathbb{B}\right)}\\[2mm]
   &\leq \beta(p-1)
   (\|u_1\|_{L^{\frac{1}{p-2}}_{(p-2)n}\left(\mathbb{B}\right)}
   +\|u_2\|_{L^{\frac{1}{p-2}}_{(p-2)n}\left(\mathbb{B}\right)})^{p-2}
    \|{|\nabla_{\mathbb{B}}(u_1-u_2)}\|_{L^{\frac{n}{2}}_{2}\left(\mathbb{B}\right)}\\[2mm]
   &\leq C_{\mathcal{B}}\|\nabla_{\mathbb{B}}(u_1-u_2)\|_{L^{\frac{n}{2}}_{2}\left(\mathbb{B}\right)}\leq C_{\mathcal{B}}\|U_1-U_2\|_{\mathcal{H}}.
  \end{split}
\end{equation}
Since  $\|\mathbb{F}(U_1)-\mathbb{F}(U_2)\|_{\mathcal{H}}=\|f(u_1)-f(u_2)\|_{L^{\frac{n}{2}}_{2}\left(\mathbb{B}\right)}$, it follows from \eqref{3.10} that $\mathbb{F}$ is locally Lispschitz in $\mathcal{H}$.

In light of the above results, by \cite[Theorem 2.5.5]{zheng2004nonlinear} or \cite[Theorem 7.2]{chue2002cpde}, we can obtain: when $U_0=(u_0, u_1)\in \mathcal{H}$, the abstract Cauchy problem \eqref{3.4} admits a unique local mild solution $U\in C \left( [0, T_{max}); \mathcal{H}\right)$ which is given by
$$U(t)=e^{\mathbb{A}t}U_0+\int_0^t e^{\mathbb{A}(t-s)}\mathbb{F}(U(s))ds,$$
defined in a maximal interval $[0, T_{max})$. When $U_0\in \mathcal{D}(\mathbb{A})$, the problem \eqref{3.4} admits a unique
strong solution $U\in C \left( [0, T_{max}); \mathcal{D}(\mathbb{A})\right)$. Moreover, if   $T_{max}<\infty$, we have $\lim_{t\to T_{max}^{-}}\|U(t)\|_{\mathcal{H}}=\infty. $

On the other hand, given $T > 0$ and any $t \in (0, T)$, we consider two mild solutions
$U^1$ and $U^2$ with initial data $U^1(0)$ and $U^2(0)$, respectively. By \eqref{3.5}, we have
\begin{equation}
    \begin{aligned}
    \|U^{1}(t)-U^{2}(t)\|_{\mathcal{H}}&\leq\quad\|e^{t\mathbb{A}}(U^{1}(0)-U^{2}(0))\|_{\mathcal{H}}
    \\&+\int_0^t\|e^{(t-s)\mathbb{A}}(\mathbb{F}(U^1(s))-\mathbb{F}(U^2(s)))\|_{\mathcal{H}}ds.
    \end{aligned}
\end{equation}
Using the local Lipschtiz property of $\mathbb{F}$, we get
$$\|U^1(t)-U^2(t)\|_{\mathcal H}\leq\|U^1(0)-U^2(0)\|_{\mathcal H}+C_0\int_0^t\|U^1(s)-U^2(s)\|_{\mathcal H}ds.$$
Then applying  Gronwall's inequality, we can obtain for any $t \in [0, T]$,
$$\|U^1(t)-U^2(t)\|_\mathcal{H}\leq e^{C_0T}\|U^1(0)-U^2(0)\|_\mathcal{H},$$
where $C_0=C(U^1(0),U^2(0))$. Thus, we deduce the continuous dependence of the mild solution on the initial data. The proof is complete.
\end{proof}

\section{ Global existence and energy decay}
First, we prove the  local solution obtained in Theorem \ref{zz} is global  provided the damping term  dominates the source, i.e., $m\leq p$ by the continuation principle.  More precisely,
\begin{theorem} (Global existence for $p\leq m$ )\\
Suppose that $m\geq p>2$, $u_{0}\in H_{2,0}^{1,\frac{n}{2}}(\mathbb{B})$ and $u_{1}\in {L^{\frac{n}{2}}_{2}}(\mathbb{B})$. Let $u(t)$ be the local solution to problem \eqref{1.1} on $[0,T_{max})$, then $u(t)$ is a global solution to problem \eqref{1.1}.
\end{theorem}
\begin{proof}
In order to show that the local solution obtained in Theorem \ref{zz} is global, by  the continuation principle, we only need to  show that $\|u_{t}\|^{2}_{L^{\frac{n}{2}}_{2}\left(\mathbb{B}\right)}
+\|\nabla_{\mathbb{B}}u\|^{2}_{L_{2}^{\frac{n}{2}}\left(\mathbb{B}\right)}$  is uniformly bounded with respect to $t$ on any interval $[0,T]$. To this end, we introduce the following two Lyapunov functionals  $$e(t)=\frac{1}{2}\|u_t\|^2_{L^{\frac{n}{2}}_{2}\left(\mathbb{B}\right)}+
\frac{1}{2}\|\nabla_{\mathbb{B}}u_{t}\|^{2}_{L_{2}^{\frac{n}{2}}\left(\mathbb{B}\right)}
-\frac{1}{2}\gamma\|V(x)^{\frac{1}{2}}u\|^{2}_{L_{2}^{\frac{n}{2}}\left(\mathbb{B}\right)},$$
$$E_{1}(t)=e(t)+\frac{1}{p}\|g(x)^{\frac{1}{p}}u\|^{2}_{L_{p}^{\frac{n}{p}}\left(\mathbb{B}\right)}.$$
From the energy identity \eqref{2.7}, we get that $$e(t)+\int_{0}^{t}(|u_\tau|^{m-2}u_\tau,u_\tau)d\tau=e(0)+\int_{0}^{t}(g(x)|u|^{p-2}u,u_\tau)d\tau.$$
Thus, $$e'(t)=-\|u_t\|^m_{L^{\frac{n}{m}}_{m}\left(\mathbb{B}\right)}+
\int_{\mathbb{B}}g(x)|u|^{p-1}u_{t}\frac{dx_{1}}{x_{1}}dx'.$$
Then, we have $$E'_{1}(t)=-\|u_t\|^m_{L^{\frac{n}{m}}_{m}\left(\mathbb{B}\right)}+
2\int_{\mathbb{B}}g(x)|u|^{p-1}u_{t}\frac{dx_{1}}{x_{1}}dx'$$
Now, we estimate the last term by cone Sobolev embedding and cone H\"{o}lder inequality
\begin{equation*}
 \begin{split}&2\int_{\mathbb{B}}g(x)|u|^{p-1}u_{t}\frac{dx_{1}}{x_{1}}dx'\\[2mm]
   &\leq 2\int_{\mathbb{B}}|g(x)|u|^{p-1}u_t|\frac{dx_{1}}{x_{1}}dx'\\[2mm]
   &\leq 2\|g(x)^\frac{1}{p}u\|^{p-1}_{L_{p}^{\frac{n}{p}}\left(\mathbb{B}\right)}
   {\|g(x)^\frac{1}{p}u_{t}\|_{L_{p}^{\frac{n}{p}}\left(\mathbb{B}\right)}}\\[2mm]
   &\leq C\|g(x)^\frac{1}{p}u\|^{p-1}_{L_{p}^{\frac{n}{p}}\left(\mathbb{B}\right)}
   {\|u_{t}\|_{L_{p}^{\frac{n}{p}}\left(\mathbb{B}\right)}}\\[2mm]
   &\leq C\|g(x)^\frac{1}{p}u\|^{p-1}_{L_{p}^{\frac{n}{p}}\left(\mathbb{B}\right)}
   {\|u_{t}\|_{L_{m}^{\frac{n}{m}}\left(\mathbb{B}\right)}},\\[2mm]
    \end{split}
\end{equation*}
where $C$ denote different positive constant depending on the known constants that may have different values in different places.
It follows from Young's inequality and $\frac{(p-1)m}{m-1}\leq p$ ( since $p\leq m$) that
\begin{equation*}
   \begin{split}&C\|g(x)^\frac{1}{p}u\|_{L_{p}^{\frac{n}{p}}\left(\mathbb{B}\right)}^{p-1}
   {\|u_{t}\|_{L_{m}^{\frac{n}{m}}\left(\mathbb{B}\right)}}\\[2mm]
  &\leq C\|g(x)^\frac{1}{p}u\|_{L_{p}^{\frac{n}{p}}\left(\mathbb{B}\right)}^{\frac{(p-1)m}{m-1}}
  +\|u_{t}\|_{L_{m}^{\frac{n}{m}}\left(\mathbb{B}\right)}^{m}\\[2mm]
  &\leq C\|g(x)^\frac{1}{p}u\|_{L_{p}^{\frac{n}{p}}\left(\mathbb{B}\right)}^{p}
  +\|u_{t}\|_{L_{m}^{\frac{n}{m}}\left(\mathbb{B}\right)}^{m}+C_{1}.
    \end{split}
\end{equation*}
Thus, from the above results, we have $E'_{1}(t)\leq C E_{1}(t)+C_{1}.$
Then, by  Gronwall's inequality, we have
\begin{equation*}
  E_{1}(t)\leq (E_{1}(0)+\frac{C_{1}}{C})e^{Ct}\leq (E_{1}(0)+\frac{C_{1}}{C})e^{cT}< +\infty
\end{equation*}
for any $t \in [0,T]$.
From this estimate and the continuation principle, we complete the proof of this theorem.
\end{proof}

\begin{lemma}\label{q}
If $E(0)<d$, $u_0\in \mathcal{W} $ and $u_1\in L_{2}^{\frac{n}{2}}\left(\mathbb{B}\right)$, there exists a constant $\theta \in(0,1)$ such that
 $$I(u)\geq \theta \|\nabla_{\mathbb{B}}u\|^{2}_{L_{2}^{\frac{n}{2}}\left(\mathbb{B}\right)},\,\,\,t\in[0, T_{max}).$$
\end{lemma}
\begin{proof}
Since $E(0)<d$ and $u_0\in \mathcal{W}$, we can obtain $u(t)\in \mathcal{W}$, $t\in [0, T_{max})$ by Lemma \ref{j}-(i). From Lemma \ref{b}, we have
\begin{equation*}
   \begin{split}\|g(x)^\frac{1}{p}u\|_{L_{p}^{\frac{n}{p}}\left(\mathbb{B}\right)}^{p}&\leq
   C_\ast^{p}\|\nabla_{\mathbb{B}}u\|^{p}_{L_{2}^{\frac{n}{2}}\left(\mathbb{B}\right)}\\[2mm]
   &\leq C_\ast^{p}(\frac{2pE(0)}{(p-2)(1-\gamma(C^\ast)^{2})})^{\frac{p-2}{2}}
   \|\nabla_{\mathbb{B}}u\|^{2}_{L_{2}^{\frac{n}{2}}\left(\mathbb{B}\right)}.
   \end{split}
\end{equation*}
Therefore, we get from the definition of $I(u)$ and Lemma \ref{d} that,
\begin{equation*}
   \begin{split}I(u)&\geq
   \|\nabla_{\mathbb{B}}u\|^{2}_{L_{2}^{\frac{n}{2}}\left(\mathbb{B}\right)}
   -\gamma(C^\ast)^{2}\|\nabla_{\mathbb{B}}u\|^{2}_{L_{2}^{\frac{n}{2}}\left(\mathbb{B}\right)}
   -C_\ast^{p}(\frac{2pE(0)}{(p-2)(1-\gamma(C^\ast)^{2})})^{\frac{p-2}{2}}\|\nabla_{\mathbb{B}}u\|^{2}_{L_{2}^{\frac{n}{2}}\left(\mathbb{B}\right)}\\[2mm]
   &\geq [1-\gamma(C^\ast)^{2}-C_\ast^{p}(\frac{2pE(0)}{(p-2)(1-\gamma(C^\ast)^{2})})^{\frac{p-2}{2}}]\|\nabla_{\mathbb{B}}u\|^{2}_{L_{2}^{\frac{n}{2}}\left(\mathbb{B}\right)}.
   \end{split}
\end{equation*}
Due to $E(0)<d $, we have $\theta:=[1-\gamma(C^\ast)^{2}-C_\ast^{p}(\frac{2pE(0)}{(p-2)(1-\gamma(C^\ast)^{2})})^{\frac{p-2}{2}}]>0$.
Thus, we have that $I(u)\geq \theta \|\nabla_{\mathbb{B}}u_{t}\|^{2}_{L_{2}^{\frac{n}{2}}\left(\mathbb{B}\right)} $.
\end{proof}

\begin{theorem}\label{A} (Global existence and energy decay for $E(0)<d$ )\\
Let $\gamma\in[0,\frac{1}{(C^{\ast})^{2}})$,  $u_{0}\in H_{2,0}^{1,\frac{n}{2}}(\mathbb{B})$,  $u_{1}\in {L^{\frac{n}{2}}_{2}}\left(\mathbb{B}\right)$ and $u(t)$ be the  local solution to problem \eqref{1.1} on $[0,T_{max})$. If $ I(u_0)\geq0$ and $E(0)<d$, then $T_{max}=\infty$, i.e., $u(t)$ exists globally and  $u(t)\in \mathcal{W}$ for all $t\in[0, \infty).$  Moreover, it has the following
energy decay estimate:
\begin{equation}\label{4.1}
\begin{split}
\qquad &E(t)\leq Ke^{-\kappa t},\,\,\,\text{for}\,\,m=2;\\[2mm]
  E(t)&\leq \left(E(0)^\frac{m-2}{m}+\frac{(m-2)\tau}{2 }[t-1]^{+}\right)^\frac{-2}{m-2},\,\,\,\text{for},\,\,m>2,
\end{split}\end{equation}
 where the positive constants $K$, $\kappa$  and $\tau$
are  given in the proof.\end{theorem}

\begin{proof}
We divide the proof into two steps: Firstly, we show the global existence of solutions. Secondly, we prove the energy decay estimate.

\textbf{Step 1: Global existence.}
It suffices to show that $\|u_{t}\|^{2}_{L^{\frac{n}{2}}_{2}\left(\mathbb{B}\right)}
+\|\nabla_{\mathbb{B}}u\|^{2}_{L_{2}^{\frac{n}{2}}\left(\mathbb{B}\right)}<+\infty$ with respect to $t$. We splits the proof into two cases.

Case 1: $I(u_0)>0$ and $E(0)<d$. It is easy to get that $u_0\in \mathcal{W}$ by \eqref{2.5} and the definition of $\mathcal{W}$.   Then, it follows from Lemma \ref{j}-(i) that $u(t)\in \mathcal{W}$ for any $t\in [0, T_{max}).$  Therefore, combining the \eqref{2.1},\eqref{2.3} and Lemma \ref{e}, we obtain
\begin{equation*}
  \begin{split} d> E(0)&= E(t)+\int_{0}^{t}\|u_{\tau}\|^{m}_{L^{\frac{n}{m}}_{m}\left(\mathbb{B}\right)}d\tau\\[2mm]
  &=\frac{1}{2}\|u_t\|^2_{L^{\frac{n}{2}}_{2}\left(\mathbb{B}\right)}+J(u)
  +\int_{0}^{t}\|u_{\tau}\|^{m}_{L^{\frac{n}{m}}_{m}\left(\mathbb{B}\right)}d\tau\\[2mm]
  &=\frac{1}{2}\|u_t\|^2_{L^{\frac{n}{2}}_{2}\left(\mathbb{B}\right)}+\frac{1}{p}I(u)
  +\frac{p-2}{2p}
  (\|\nabla_{\mathbb{B}}u\|^{2}_{L_{2}^{\frac{n}{2}}\left(\mathbb{B}\right)}
  -\gamma\|V(x)^{\frac{1}{2}}u\|^{2}_{L_{2}^{\frac{n}{2}}\left(\mathbb{B}\right)})\\[2mm]
  &\geq\frac{1}{2}\|u_t\|^2_{L^{\frac{n}{2}}_{2}\left(\mathbb{B}\right)}+\frac{1}{p}I(u)
  +\frac{p-2}{2p}(1-\gamma(C^\ast)^{2})
  \|\nabla_{\mathbb{B}}u\|^{2}_{L_{2}^{\frac{n}{2}}\left(\mathbb{B}\right)}\\[2mm]
  &\geq\frac{1}{2}\|u_t\|^2_{L^{\frac{n}{2}}_{2}\left(\mathbb{B}\right)}
  +\frac{p-2}{2p}(1-\gamma(C^\ast)^{2})\|\nabla_{\mathbb{B}}u\|^{2}
  _{L_{2}^{\frac{n}{2}}\left(\mathbb{B}\right)},
  \end{split}
\end{equation*}
which yields that $\|u_{t}\|^{2}_{L^{\frac{n}{2}}_{2}\left(\mathbb{B}\right)}
+\|\nabla_{\mathbb{B}}u\|^{2}_{L_{2}^{\frac{n}{2}}\left(\mathbb{B}\right)}<\infty$ with respect to $t$. In view of Theorem \ref{zz},  by the continuation principle, the above estimate yields that $T_{max}=+\infty$, i.e.  the solution $u(t)$ exists globally.

Case 2: $I(u_0)=0$ and $E(0)<d$. We claim  that $u_0=0$. Indeed, if $u_0\neq0$, then $u_0\in \mathcal{N}$, which implies that $J(u_0)\geq d$. Then it follows from \eqref{2.5} that $E(0)\geq d$, which is a contradiction with $E(0)<d$.   Thus $u_0=0\in\mathcal{W}$. The remaining proof  is the same as case 1.

 \textbf{Step 2: Energy decay.}
 By integrating $\frac{d}{dt}E(t)=-\|u_t\|^m_{L^{\frac{n}{m}}_{m}\left(\mathbb{B}\right)}$ over $[t,t+1]$ for $t\geq0$, we have that
 \begin{equation}\label{4.2}
   E(t)-E(t+1)=\int_{t}^{t+1}\|u_{\tau}\|^{m}_{L^{\frac{n}{m}}_{m}\left(\mathbb{B}\right)}d\tau:=D^{m}(t).
\end{equation}
 By  cone H\"{o}lder inequality,   we deduce
 \begin{equation*} \label{eq4.9}
   \begin{split}\int_{t}^{t+1}\|u_{\tau}\|^{2}_{L^{\frac{n}{2}}_{2}\left(\mathbb{B}\right)}d\tau
   &\leq C \int_{t}^{t+1}\|u_{\tau}\|^{2}_{L^{\frac{n}{m}}_{m}\left(\mathbb{B}\right)}d\tau\\[2mm]
   &\leq C\left( {\int_{t}^{t+1}{1}^{\frac{m}{m-2}}d\tau}\right)^{\frac{m-2}{m}}
   \left({\int_{t}^{t+1}\|u_{\tau}\|^{m}_{L^{\frac{n}{m}}_{m}\left(\mathbb{B}\right)}d\tau}\right)^{\frac{2}{m}}\\[2mm]
   &= C D^{2}{(t)}.
   \end{split}
\end{equation*}
By the mean value theorem, we deduce that there exists $t_{1}\in[t,t+\frac{1}{4}]$, $t_{2}\in[t+\frac{3}{4}, t+1]$,
such that$${\|u_t(t_i)\|^2_{L^{\frac{n}{2}}_{2}(\mathbb{B})}}\leq 4 C D^{2}{(t)}.$$
Now, multiplying \eqref{1.1} by $u$, integrating it over $\mathbb{B} \times[t_{1}, t_{2}]$ and integrating by parts, we deduce that
$$\int_{t_{1}}^{t_{2}}I(u)dt=-\int_{t_{1}}^{t_{2}}\int_{\mathbb{B}}u_{tt}u\frac{dx_{1}}{x_{1}}dx'dt
-\int_{t_{1}}^{t_{2}}\int_{\mathbb{B}}|u_{t}|^{m-2}u_{t}u\frac{dx_{1}}{x_{1}}dx'dt.$$
Then, we estimate the right side of the above equation,
\begin{equation} \label{eq4.10}
   \begin{split}&\int_{t_{1}}^{t_{2}}\int_{\mathbb{B}}u_{tt}u\frac{dx_{1}}{x_{1}}dx'dt\\[2mm]
   &= \int_{\mathbb{B}}u_{t}u\mid^{t_{2}}_{t_{1}}\frac{dx_{1}}{x_{1}}dx'
   -\int_{\mathbb{B}}\int_{t_{1}}^{t_{2}}|u_{t}|^{2}dt\frac{dx_{1}}{x_{1}}dx'\\[2mm]
   &\leq\mathop{\sum}\limits_{i=1}^{2}\|u_t(t_{i})\|_{L^{\frac{n}{2}}_{2}\left(\mathbb{B}\right)}\|u(t_{i})\|_{L^{\frac{n}{2}}_{2}\left(\mathbb{B}\right)}
   +\int_{t_{1}}^{t_{2}}\|u_t\|^2_{L^{\frac{n}{2}}_{2}\left(\mathbb{B}\right)}dt.
   \end{split}
\end{equation}
From the Lemma \ref{a}, ${L^{\frac{n}{m}}_{m}\left(\mathbb{B}\right)}\hookrightarrow{L^{\frac{n}{2}}_{2}\left(\mathbb{B}\right)}$ and
$\lambda_1^2\|u\|^2_{L^{\frac{n}{2}}_{2}\left(\mathbb{B}\right)}\leq \|\nabla_\mathbb{B}u\|^2_{L^{\frac{n}{2}}_{2}\left(\mathbb{B}\right)}$, we have
\begin{equation} \label{eq4.11}
   \begin{split}&\|u(t_i)\|_{L^{\frac{n}{2}}_{2}\left(\mathbb{B}\right)}\|u_{t}(t_i)\|
   _{L^{\frac{n}{2}}_{2}\left(\mathbb{B}\right)}\\[2mm]
   &\leq 2CD(t)\frac{1}{\lambda_1}\|\nabla_\mathbb{B}u(t_{i})\|_{L^{\frac{n}{2}}_{2}\left(\mathbb{B}\right)}\\[2mm]
   &\leq 2CD(t)\frac{1}{\lambda_1}
   \left(\frac{2p}{(p-2)[1-\gamma(C^\ast)^{2}]}\right)^{\frac{1}{2}}\mathop{sup}\limits_{t_{1}\leq s\leq t_{2}}E^{\frac{1}{2}}(s),
   \end{split}
\end{equation}
and
\begin{equation} \label{eq4.12}
   \begin{split}\int_{t_{1}}^{t_{2}}\|u_t\|^2_{L^{\frac{n}{2}}_{2}\left(\mathbb{B}\right)}dt
   &\leq C\int_{t_{1}}^{t_{2}}\|u_t\|^2_{L^{\frac{n}{m}}_{m}\left(\mathbb{B}\right)}dt\\[2mm]
   &\leq C\left({\int_{t_{1}}^{t_{2}}\|u_t\|^m_{L^{\frac{n}{m}}_{m}\left(\mathbb{B}\right)}}
   dt\right)^\frac{2}{m}\\[2mm]
   &\leq CD^{2}(t).
   \end{split}
\end{equation}
By  $2<m<p$ and cone Sobolev imbedding, we observe that
\begin{equation} \label{eq4.13}
   \begin{split}&|-\int_{t_{1}}^{t_{2}}\int_{\mathbb{B}}|u_t|^{m-2}u_t u|\frac{dx_{1}}{x_{1}}dx'dt|\\[2mm]
   &\leq \int_{t_{1}}^{t_{2}}\int_{\mathbb{B}}|u_t|^{m-1} u\frac{dx_{1}}{x_{1}}dx'dt\\[2mm]
   &\leq \int_{t_{1}}^{t_{2}}
   \int_{\mathbb{B}}\left(|u_t|^{m}
   \frac{dx_{1}}{x_{1}}dx'\right)^\frac{m-1}{m}
   \left(\int_{\mathbb{B}}|u|^{m}\frac{dx_{1}}{x_{1}}dx'\right)^\frac{1}{m}dt\\[2mm]
   &\leq C\int_{t_{1}}^{t_{2}}\|u_t\|^{m-1}_{L^{\frac{n}{m}}_{m}\left(\mathbb{B}\right)}\|\nabla_{\mathbb{B}}u\|
   _{L_{2}^{\frac{n}{2}}\left(\mathbb{B}\right)}dt\\[2mm]
   &\leq C D(t)
   \left(\frac{2p}{(p-2)[1-\gamma(C^\ast)^{2}]}\right)^{\frac{1}{2}}\mathop{sup}\limits_{t_{1}\leq s\leq t_{2}}E^{\frac{1}{2}}(s)\int_{t_1}^{t_2}\|u_{\tau}\|^{m-1}_{L^{\frac{n}{m}}_{m}\left(\mathbb{B}\right)
}d\tau\\[2mm]
   &\leq C D(t)
   \left(\frac{2p}{(p-2)[1-\gamma(C^\ast)^{2}]}\right)^{\frac{1}{2}}\mathop{sup}\limits_{t_{1}\leq s\leq t_{2}}E^{\frac{1}{2}}(s)D^{m-1}(t).
   \end{split}
\end{equation}
Hence, $\int_{t_{1}}^{t_{2}}I(u)dt$ can be estimated  as
\begin{equation*}
\begin{split}\int_{t_{1}}^{t_{2}}I(u)dt
&\leq 4CD(t)\frac{1}{\lambda_1}
   \left(\frac{2p}{(p-2)[1-\gamma(C^\ast)^{2}]}\right)^{\frac{1}{2}}\mathop{sup}\limits_{t_{1}\leq s\leq t_{2}}E^{\frac{1}{2}}(s)+CD^2(t)\\[2mm]
   &+\left(\frac{2p}{(p-2)[1-\gamma(C^\ast)^{2}]}\right)^{\frac{1}{2}}\mathop{sup}\limits_{t_{1}\leq s\leq t_{2}}E^{\frac{1}{2}}(s)D^{m-1}(t).
\end{split}
\end{equation*}
Meanwhile, by Lemma \ref{q} we deduce that
\begin{equation} \label{eq4.14}
   \begin{split}J(u)&=(\frac{1}{2}-\frac{1}{p})
   \left(\|\nabla_{\mathbb{B}}u\|^{2}_{L_{2}^{\frac{n}{2}}\left(\mathbb{B}\right)}-
   \gamma\|V(x)^{\frac{1}{2}}u\|^{2}_{L_{2}^{\frac{n}{2}}\left(\mathbb{B}\right)}\right)
   +\frac{1}{p}I(u)\\[2mm]
   &\leq (\frac{1}{2}-\frac{1}{p})
   \|\nabla_{\mathbb{B}}u\|^{2}_{L_{2}^{\frac{n}{2}}\left(\mathbb{B}\right)}
   +\frac{1}{p}I(u)\\[2mm]
   &\leq \left(\frac{1}{p}+\frac{p-2}{2p\theta}\right)I(u).
   \end{split}
\end{equation}
Hence, we have
\begin{equation} \label{eq4.15}
   \begin{split}E(t)&=
   \frac{1}{2}\|u_t\|^2_{L^{\frac{n}{2}}_{2}\left(\mathbb{B}\right)}
   +J(u)\\[2mm]
   &\leq\frac{1}{2}\|u_t\|^2_{L^{\frac{n}{2}}_{2}\left(\mathbb{B}\right)}
   +C_2 I(u)
   \end{split}
\end{equation}
where $C_{2}=\left(\frac{1}{p}+\frac{p-2}{2p\theta}\right)$.
Integrating $E(t)=\frac{1}{2}\|u_t\|^2_{L^{\frac{n}{2}}_{2}\left(\mathbb{B}\right)}+J(u)$ over $(t_1,t_2)$, we have
\begin{equation} \label{eq4.16}
   \begin{split}\int_{t_{1}}^{t_{2}}E(t)dt&\leq\int_{t_{1}}^{t_{2}}\frac{1}{2}\|u_t\|^2_{L^{\frac{n}{2}}_{2}\left(\mathbb{B}\right)}dt
   +C_{2}\int_{t_{1}}^{t_{2}}I(u)dt\\[2mm]
   &\leq \frac{1}{2}D^{2}(t)+C_{2}
   [4CD(t)\frac{1}{\lambda_1}
   \left(\frac{2p}{(p-2)[1-\gamma(C^\ast)^{2}]}\right)^{\frac{1}{2}}\mathop{sup}\limits_{t_{1}\leq s\leq t_{2}}E^{\frac{1}{2}}(s)]\\[2mm]
   &+C_2CD^2(t)+
   C_2\left(\frac{2p}{(p-2)[1-\gamma(C^\ast)^{2}]}\right)^{\frac{1}{2}}\mathop{sup}\limits_{t_{1}\leq s\leq t_{2}}E^{\frac{1}{2}}(s)D^{m-1}(t).
   \end{split}
\end{equation}
On the other hand,  integrating  $E'(t)=-\|u_{t}\|^{m}_{L^{\frac{n}{m}}_{m}}$ over $[t,t_{2})$, we obtain, $$E(t)=E(t_{2})+\int_{t}^{t_2}\|u_{\tau}\|^{m}_{L^{\frac{n}{m}}_{m}}d\tau.$$
Noticing $t_{2}-t_{1}\geq\frac{1}{2}$, we see that $$2\int_{t_{1}}^{t_{2}}E(t)d\tau\geq E(t_{2}).$$  Then, combining the definition of $D^m(t)$ and \eqref{eq4.16}, we see that
\begin{equation*}
   \begin{split}E(t)&=E(t_{2})+\int_{t}^{t_2}\|u_{\tau}\|^{m}_{L^{\frac{n}{m}}_{m}\left(\mathbb{B}\right)}d\tau\\[2mm]
   &\leq 2\int_{t_{1}}^{t_{2}}E(t)dt+\int_{t_{1}}^{t_{2}}\|u_{\tau}\|^{m}_{L^{\frac{n}{m}}_{m}\left(\mathbb{B}\right)}d\tau\\[2mm]
   &\leq 2\int_{t_{1}}^{t_{2}}E(t)dt+D^{m}(t)\\[2mm]
   &\leq D^{2}(t)+2C_2CD^2(t)+D^{m}(t)\\[2mm]
   &+C_{2}[8CD(t)\frac{1}{\lambda_1}
   +2C_2D^{m-1}]\left(\frac{2p}{(p-2)[1-\gamma(C^\ast)^{2}]}\right)^{\frac{1}{2}}
   \mathop{sup}\limits_{t_{1}\leq s\leq t_{2}}E^{\frac{1}{2}}(s)\\[2mm]
   &= \left(1+2 C C_2\right)D^2(t)+C_2(D(t)+D^{m-1}(t))\mathop{sup}\limits_{t_{1}\leq s\leq t_{2}}E^{\frac{1}{2}}(s)+D^{m}(t).
   \end{split}
\end{equation*}
Hence, it follows from Young's inequality that
\begin{equation}\label{4.9}
 \begin{split}E(t)&\leq C_{3}[D^{2}(t)+D^{2(m-1)}(t)+D^{m}(t)] \\[2mm]
 &\leq C_{3}[1+D^{2(m-2)}(t)+D^{m-2}(t)]D^{2}(t),
 \end{split}
\end{equation}
 where $C_{3}>1$. Then, it follows from \eqref{4.2} and \eqref{4.9} that
$$E(t)\leq C_{3}[1+E^{\frac{2(m-1)}{m}}(0)+E^{\frac{(m-2)}{m}}(0)]D^{2}(t),$$
which implies that $$E^{\frac{m}{2}}(t)\leq C_{4}^{\frac{m}{2}}D^{m}(t)\leq C_{4}^{\frac{m}{2}}(E(t)-E(t+1)),$$
where $C_{4}=C_4(E(0)):=C_{3}[1+E^{\frac{2(m-1)}{m}}(0)+E^{\frac{(m-2)}{m}}(0)]>1,$ $\displaystyle \lim_{E(0)\to 0^+}C_4(E(0))=C_3.$
Then, by Lemma \ref{m} , we get the energy decay estimate  \eqref{4.1} holds $\tau= (C_4(E(0)))^{-\frac{m}{2}}$.
\end{proof}

\begin{theorem} \label{B}(Global existence and energy decay for $E(0)=d$ )\\
Let $\gamma\in[0,\frac{1}{(C^{\ast})^{2}})$, $u_{0}\in H_{2,0}^{1,\frac{n}{2}}(\mathbb{B})$,  $u_1\in L^{\frac{n}{2}}_{2}\left(\mathbb{B}\right)$ and $u(t)$ be the  local solution to problem \eqref{1.1} on $[0,T_{max})$. If  $I(u_0)\geq0$, $\|u_1\|_{L_2^{\frac{n}{2}}(\mathbb{B})}>0$ and $E(0)=d$, then $T_{max}=\infty$, i.e., $u(t)$ exists globally.  Moreover, it has the following
energy decay estimate:
\begin{equation*}
\begin{split}\qquad &E(t)\leq K'E^{-\kappa't},\,\,\,\text{for}\,\,m=2;\\[2mm]
  E(t)&\leq \left(E(0)^\frac{m-2}{m}+\frac{(m-2)\tau'}{2 }[t-1]^{+}\right)^\frac{-2}{m-2},\,\,t\geq t_1,\,\,\,\text{for}\,\,m>2,
\end{split}\end{equation*}
 for some $t_1>0$ and some positive constants $K', \kappa'$ and $\tau'$.
\end{theorem}
\begin{proof} We split the condition $I(u_0)\geq0$ into two cases.

Case 1: $I(u_0)=0$. Since  $\|u_1\|_{L_2^{\frac{n}{2}}(\mathbb{B})}>0$,  by \eqref{2.5} and $E(0)=d$, we can see that $J(u_0)<d.$ Then we have $u_0=0\in \mathcal{W}$ by the similar argument as case 2 in the proof  of Theorem \ref{A}. The remaining proof is the same as Theorem \ref{A}.

Case 2: $I(u_0)>0$. We claim that there a sufficiently  small time $t_0>0$ such that $E(t)<d$ for all $t\in (0, t_0]$. Indeed, it the claim is not true, noticing that $E(t)$ is non-increasing with respect to $t$, there must exist a decreasing sequence $\{t_k\}_{k=1}^{\infty}$ ($t_k\to 0$ as $k\to\infty$) such that $E(t_k)=d$.  Then, it follows from \eqref{2.7} that $\int_0^{t_k}\|u_t\|^m_{L^{\frac{n}{m}}_{m}\left(\mathbb{B}\right)}dt=0,\,\,k=1,2,\cdots.$ Since $u_t\in C\left([0, T_{max}, L^{\frac{n}{2}}_{2}\left(\mathbb{B}\right)\right)$, by cone H\"{o}lder inequality, we see that $\|u_t(t_k)\|_{L^{\frac{n}{2}}_{2}\left(\mathbb{B}\right)}=0,\,\, k=1,2,\cdots.$  Then, letting $k\to\infty$, we get $\|u_1\|_{L^{\frac{n}{2}}_{2}\left(\mathbb{B}\right)}=0$, which contradicts with the condition $\|u_1\|_{L^{\frac{n}{2}}_{2}\left(\mathbb{B}\right)}>0.$  On the other hand,  it is easy to see that there exists a sufficiently small $t_1\in (0, t_0]$ such that $I(u(t_1))>0$. It follows from \eqref{2.5} that $J(u(t_1))\leq E(t_1)<d$. The above analysis shows that $u(t_1)\in \mathcal{W}.$
Similarly as the proof of Theorem \ref{A}, we derive the conclusion.
\end{proof}

\section{Blow up }
In this section, we will establish that the  solution of problem \eqref{1.1} blows up at a
finite time under different initial energy levels.
\begin{theorem}\label{th5.1} (Blow up for  $E(0)<0$ )\\
Let $\gamma\in[0,\frac{1}{(C^{\ast})^{2}})$, $u_{0}\in H_{2,0}^{1,\frac{n}{2}}(\mathbb{B})$,  $u_1\in L^{\frac{n}{2}}_{2}\left(\mathbb{B}\right)$ and $u(t)$ be the  local solution to problem \eqref{1.1} on $[0,T_{max})$. If the initial energy   $E(0)<0$, then the solution to  problem \eqref{1.1} blows up at a  finite time.
\end{theorem}
\begin{proof}
Since $E(0)<0$, we derive from Lemma \ref{j} that $u(t)\in \mathcal{V}$ for all $t\in [0, T_{max})$. Suppose $u(t)$ is the global solution of \eqref{1.1}. For any given  $T>0$, setting $F(t)=\|u(t)\|^2_{L^{\frac{n}{2}}_{2}\left(\mathbb{B}\right)}$ for $t\in [0,T]$,  we have  $F'(t)=2(u,u_t)$ and
$\frac{d}{dt}(u, u_t)=\|u_t\|^2_{L^{\frac{n}{2}}_{2}\left(\mathbb{B}\right)}+\langle u_{tt}, u\rangle$.
Thus, we have$$(u,u_t)=(u_t(0),u(0))+\int_{0}^{t}(\|u_t\|^2_{L^{\frac{n}{2}}_{2}\left(\mathbb{B}\right)}+\langle u_{tt},u\rangle)d\tau$$
and
\begin{equation*}
   \begin{split}F''(t)&=2\frac{d}{dt}(u,u_t)\\[2mm]
   &= 2(\|u_t\|^2_{L^{\frac{n}{2}}_{2}\left(\mathbb{B}\right)}+\langle \Delta_\mathbb{B}u+\gamma V(x)u-|u_t|^{m-2}u_t+g(x)|u|^{p-2}u ,u \rangle)\\[2mm]
   &=4\|u_t\|^2_{L^{\frac{n}{2}}_{2}\left(\mathbb{B}\right)}
   -4E(t)+\frac{2(p-2)}{p}\|g(x)^\frac{1}{p}u\|_{L_{p}^{\frac{n}{p}}\left(\mathbb{B}\right)}^{p}
   -2\int_{\mathbb{B}}|u_t|^{m-2}u_t u\frac{dx_{1}}{x_{1}}dx'.
   \end{split}
\end{equation*}

Let $H(t)=-E(t)$, from the energy equation \eqref{2.7}, we deduce that $$H'(t)=\|u_{t}\|^{m}_{L^{\frac{n}{m}}_{m}}\geq0,$$
which implies
\begin{equation}\label{5.1}
  0<H(0)\leq H(t)\leq \frac{1}{p}\|g(x)^\frac{1}{p}u\|_{L_{p}^{\frac{n}{p}}\left(\mathbb{B}\right)}^{p},\,\,t\in [0,T].
\end{equation}
Now, let $L(t)=H^{1-\eta}(t)+\varepsilon F'(t)$, where $\varepsilon>0$, $\eta>0$, we deduce that
\begin{equation*}
  \begin{split}
  L'(t)=&(1-\eta)H^{-\eta}(t)H'(t)+4\varepsilon\|u_t\|^2_{L^{\frac{n}{2}}_{2}\left(\mathbb{B}\right)}+4\varepsilon H(t)
+\frac{2\varepsilon(p-2)}{p}\|g(x)^\frac{1}{p}u\|_{L_{p}^{\frac{n}{p}}\left(\mathbb{B}\right)}^{p}\\[2mm]
&-2\varepsilon\int_{\mathbb{B}}|u_t|^{m-2}u_t u\frac{dx_{1}}{x_{1}}dx'.
  \end{split}
\end{equation*}
By Young's inequality, we have
\begin{equation*}
   \begin{split}|\int_{\mathbb{B}}|u_t|^{m-2}u_t u\frac{dx_{1}}{x_{1}}dx'|&\leq\|u_t\|^{m-1}_{L^{\frac{n}{m}}_{m}\left(\mathbb{B}\right)}\|u\|_{L^{\frac{n}{m}}_{m}\left(\mathbb{B}\right)}\\[2mm]
   &\leq \frac{m-1}{m}\delta^{\frac{-m}{m-1}}\|u_t\|^m_{L^{\frac{n}{m}}_{m}\left(\mathbb{B}\right)}
   +\frac{1}{m}\delta^{m}\|u\|^m_{L^{\frac{n}{m}}_{m}\left(\mathbb{B}\right)}.
   \end{split}
\end{equation*}
Thus, we arrive at
\begin{equation} \label{5.2}
   \begin{split}L'(t)&\geq(1-\eta)H^{-\eta}(t)H'(t)+4\varepsilon\|u_t\|^2_{L^{\frac{n}{2}}_{2}\left(\mathbb{B}\right)}
   +4\varepsilon H(t)
+\frac{2\varepsilon(p-2)}{p}\|g(x)^\frac{1}{p}u\|_{L_{p}^{\frac{n}{p}}\left(\mathbb{B}\right)}^{p}\\[2mm]
&-2\varepsilon\frac{m-1}{m}\delta^{\frac{-m}{m-1}}\|u_t\|^m_{L^{\frac{n}{m}}_{m}\left(\mathbb{B}\right)}
   -2\varepsilon\frac{1}{m}\delta^{m}\|u\|^m_{L^{\frac{n}{m}}_{m}\left(\mathbb{B}\right)}.
   \end{split}
\end{equation}
Therefore, taking $\delta$ such that $\delta^{\frac{-m}{m-1}}=kH^{-\eta}(t)$ for sufficiently large $k$ to be specified later and combining ${L^{\frac{n}{p}}_{p}\left(\mathbb{B}\right)}\hookrightarrow{L^{\frac{n}{m}}_{m}\left(\mathbb{B}\right)}$
, we have
\begin{equation*}
   \begin{split}\delta^{m}\|u\|^m_{L^{\frac{n}{m}}_{m}\left(\mathbb{B}\right)}
   &=k^{-(m-1)}H^{\eta(m-1)}(t)\|u\|^m_{L^{\frac{n}{m}}_{m}\left(\mathbb{B}\right)}\\[2mm]
   &\leq k^{1-m}C(\frac{\beta}{p})^{\eta(m-1)}\|u\|^{p\eta(m-1)}_{L^{\frac{n}{p}}_{p}\left(\mathbb{B}\right)}\|u\|^m_{L^{\frac{n}{p}}_{p}\left(\mathbb{B}\right)}\\[2mm]
   &\leq k^{1-m}C(\frac{\beta}{p})^{\eta(m-1)}\|u\|^{p\eta(m-1)+m}_{L^{\frac{n}{p}}_{p}\left(\mathbb{B}\right)}.
   \end{split}
\end{equation*}
Choosing  $0<\eta\leq\frac{p-m}{p(m-1)}$ such that $s:=\eta(m-1)p+m-p\leq0$, which implies $\eta<\frac{p-m}{(m-1)p}$, we have
$\|u\|^s_{L^{\frac{n}{p}}_{p}\left(\mathbb{B}\right)}\leq (\frac{p}{\beta}H(0))^\frac{s}{p}$. Thus, $L'(t)$ becomes
\begin{equation} \label{5.3}
   \begin{split} L'(t)&\geq[(1-\eta)-\frac{2\varepsilon k(m-1)}{m}]H^{-\eta}(t)H'(t)
   +4\varepsilon H(t)
   +4\varepsilon\|u_t\|^2_{L^{\frac{n}{2}}_{2}\left(\mathbb{B}\right)}\\[2mm]
   &+\frac{2\varepsilon (p-2)}{p}\|g(x)^\frac{1}{p}u\|_{L_{p}^{\frac{n}{p}}}^{p}
   -C\frac{k^{1-m}}{m}(\frac{\beta}{p})^{\eta(m-1)}(\frac{p}{\beta})^{\frac{s}{p}}H^{\frac{s}{p}}(0)\|u\|_{L_{p}^{\frac{n}{p}}\left(\mathbb{B}\right)}^{p}\\[2mm]
   &\geq 2\varepsilon [\frac{\alpha(p-2)}{p}-C\frac{k^{1-m}}{m}(\frac{\beta}{p})^{\eta(m-1)}(\frac{p}{\beta})^{\frac{s}{p}}H^{\frac{s}{p}}(0)]\|u\|_{L_{p}^{\frac{n}{p}}\left(\mathbb{B}\right)}^{p}\\[2mm]
   &+[(1-\eta)-\frac{2\varepsilon k(m-1)}{m}]H^{-\eta}(t)H'(t)+4\varepsilon H(t)
   +4\varepsilon\|u_t\|^2_{L^{\frac{n}{2}}_{2}\left(\mathbb{B}\right)}.
   \end{split}
\end{equation}
Choosing $k$ sufficiently largy such that $[\frac{\alpha(p-2)}{p}-C\frac{k^{1-m}}{m}(\frac{\beta}{p})^{\eta(m-1)}(\frac{p}{\beta})^{\frac{s}{p}}H^{\frac{s}{p}}(0)]>0$
and then choosing $\varepsilon>0$ sufficiently small such that $[(1-\eta)-\frac{2\varepsilon k(m-1)}{m}]>0$, we have
$$L(0)=H^{1-\eta}(0)+\varepsilon F'(0)=H^{1-\eta}(0)+2\varepsilon\int_{\mathbb{B}}u_0u_1\frac{dx_{1}}{x_{1}}dx'>0.$$
Thus,
\begin{equation}\label{5.4}
  L'(t)\geq\delta\varepsilon[H(t)+\|u_t\|^2_{L^{\frac{n}{2}}_{2}\left(\mathbb{B}\right)}+\|u\|_{L_{p}^{\frac{n}{p}}\left(\mathbb{B}\right)}^{p}],
\end{equation}
where $\delta$ is the minimum of the coefficients of $H(t)$, $\|u_t\|^2_{L^{\frac{n}{2}}_{2}\left(\mathbb{B}\right)}$ and $\|u\|^2_{L^{\frac{n}{p}}_{p}\left(\mathbb{B}\right)}$ in \eqref{5.3}.

On the other hand, by H\"{o}lder inequality and Young's inequality, we deduce
\begin{equation*}
   \begin{split}|\int_{\mathbb{B}}uu_t\frac{dx_{1}}{x_{1}}dx'|^{\frac{1}{1-\eta}}
   &\leq\|u_t\|^\frac{1}{1-\eta}_{L^{\frac{n}{2}}_{2}\left(\mathbb{B}\right)}\|u\|^\frac{1}{1-\eta}_{L^{\frac{n}{2}}_{2}\left(\mathbb{B}\right)}\\[2mm]
   &\leq C \|u_t\|^\frac{1}{1-\eta}_{L^{\frac{n}{2}}_{2}\left(\mathbb{B}\right)}\|u\|^\frac{1}{1-\eta}_{L^{\frac{n}{p}}_{p}\left(\mathbb{B}\right)}\\[2mm]
   &\leq C(\|u\|^\frac{2}{1-2\eta}_{L^{\frac{n}{p}}_{p}\left(\mathbb{B}\right)}+\|u_t\|^2_{L^{\frac{n}{2}}_{2}\left(\mathbb{B}\right)}).
   \end{split}
\end{equation*}
Choosing the fact  $\eta>0$ such that $\frac{2}{1-2\eta}\leq p$, which implies $0<\eta\leq \min\left\{\frac{p-2}{2p},\frac{p-m}{(m-1)p}\right\}$, we have for any $0<\frac{2}{(1-2\eta)p}\leq1$
\begin{equation*}
\begin{split}
  \|u\|^\frac{2}{1-2\eta}_{L^{\frac{n}{p}}_{p}\left(\mathbb{B}\right)}
  &\leq(1+\frac{1}{H(0)})(\|u\|^{p}_{L^{\frac{n}{p}}_{p}}+H(0))\\[2mm]
  &\leq(1+\frac{1}{H(0)})(\|u\|^{p}_{L^{\frac{n}{p}}_{p}}+H(t)).
\end{split}
\end{equation*}
And,
$$|\int_{\mathbb{B}}uu_t\frac{dx_{1}}{x_{1}}dx'|^{\frac{1}{1-\eta}}\leq C(\|u\|^{p}_{L^{\frac{n}{p}}_{p
}\left(\mathbb{B}\right)}+\|u_t\|^2_{L^{\frac{n}{2}}_{2}\left(\mathbb{B}\right)}). $$
Thus, $$L^{\frac{1}{1-\eta}}(t)\leq C[H(t)+\|u_t\|^2_{L^{\frac{n}{2}}_{2}\left(\mathbb{B}\right)}+\|u\|_{L_{p}^{\frac{n}{p}}\left(\mathbb{B}\right)}^{p}].$$ Combining \eqref{5.4}, we have
$$L'(t)\geq C L^{\frac{1}{1-\eta}}(t).$$
Thus, $$L'(t)\geq[L^{\frac{-\eta}{1-\eta}}(0)-\frac{\eta C}{1-\eta}t]^{\frac{1-\eta}{-\eta}}\,\,\text{for}\,\, t\in[0,T],$$
since $L(0)>0$, which shows that $L(t)$ becomes infinite at  a finite time
\begin{equation*}
 T\leq T^\ast=\frac{1-\eta}{C\eta}L^\frac{-\eta}{1-\eta}(0).
\end{equation*}
\end{proof}

\begin{theorem}\label{th5.2} (Blow up for   $E(0)\leq d$ )\\
Let $\gamma\in[0,\frac{1}{(C^{\ast})^{2}})$,  $u_{0}\in H_{2,0}^{1,\frac{n}{2}}(\mathbb{B})$,  $u_{1}\in {L^{\frac{n}{2}}_{2}}\left(\mathbb{B}\right)$ and $u(t)$ be the  local solution to problem \eqref{1.1} on $[0,T_{max})$. Suppose that $p>m\geq2$, if one of the following conditions is satisfied\\
\textbf{(i)} $0\leq E(0)<d$, $I(u_0)<0$;\\
\textbf{(ii)} $E(0)=d$, $I(u_0)<0$ and  $(u_0,u_1)>0$,\\
then the solution of  problem \eqref{1.1} blows up at a finite time.
\end{theorem}
\begin{proof}
This proof is based on the  proof of Theorem \ref{th5.1}. We divide this proof into two cases.

\textbf{Case 1: $E(0)<d$.}\\
 Arguing by contradiction, we assume that the maximal existence time $T_{max}=\infty$. We choose $$d_1\in (E(0), d).$$  Firstly, we modify  $H(t)=-E(t)$ in Theorem \ref{th5.1} into $H(t)=d_1-E(t)$, that is
  \begin{equation*}
    H(t)=d_1-\frac{1}{2}\|u_t\|^2_{L^{\frac{n}{2}}_{2}\left(\mathbb{B}\right)}-\frac{1}{2}\|\nabla_{\mathbb{B}}u\|^{2}_{L_{2}^{\frac{n}{2}}\left(\mathbb{B}\right)}
+\frac{1}{2}\gamma\|V(x)^{\frac{1}{2}}u\|^{2}_{L_{2}^{\frac{n}{2}}\left(\mathbb{B}\right)}
+\frac{1}{p}\|g(x)^\frac{1}{p}u\|_{L_{p}^{\frac{n}{p}}\left(\mathbb{B}\right)}^{p}.
\end{equation*}
Noticing $E(0)<d$ and $I(u_0)<0$,  we see that $u_0\in \mathcal{V}$ by \eqref{2.5}. Then, we have $u(t)\in \mathcal{V}$ for all $t\geq0$ from Lemma \eqref{j}-(ii). Hence, from Lemma \ref{k} and the definition of $d$, we can see that
\begin{equation*}
   \begin{split}
&\frac{1}{2}\|\nabla_{\mathbb{B}}u(t)\|^{2}_{L_{2}^{\frac{n}{2}}\left(\mathbb{B}\right)}
-\frac{1}{2}\gamma\|V(x)^{\frac{1}{2}}u(t)\|^{2}_{L_{2}^{\frac{n}{2}}\left(\mathbb{B}\right)}-d_1\\[2mm]
&\geq\frac{1}{2}[1-\gamma(C^\ast)^{2}]\|\nabla_{\mathbb{B}}u(t)\|^{2}_{L_{2}^{\frac{n}{2}}\left(\mathbb{B}\right)}-d_1\\[2mm]
&>\frac{1}{2}[1-\gamma(C^\ast)^{2}]^{\frac{p}{p-2}}{C_\ast}^{-\frac{2p}{p-2}}-d\\[2mm]
&=\frac{2d}{p-2}>0.
\end{split}
\end{equation*}
From \eqref{5.1}, we get
\begin{equation*}
  0<H(0)\leq H(t)\leq\frac{1}{p}\|g(x)^\frac{1}{p}u(t)\|_{L_{p}^{\frac{n}{p}}\left(\mathbb{B}\right)}^{p}\leq\frac{\beta}{p}\|u(t)\|_{L_{p}^{\frac{n}{p}}\left(\mathbb{B}\right)}^{p}.
\end{equation*}
 Now, we consider the function $$L(t)=H^{1-\eta}(t)+\varepsilon F'(t),\,\,\,t\geq0,$$ where $\varepsilon>0$ and $\eta>0$ are sufficiently small parameter to be determined.
We can easily see that
\begin{equation} \label{5.5}
   \begin{split}L'(t)&\geq(1-\eta)H^{-\eta}(t)H'(t)+3\varepsilon\|u_t(t)\|^2_{L^{\frac{n}{2}}_{2}\left(\mathbb{B}\right)}
-\varepsilon\left(\|\nabla_{\mathbb{B}}u(t)\|^{2}_{L_{2}^{\frac{n}{2}}\left(\mathbb{B}\right)}
-\gamma\|V(x)^{\frac{1}{2}}u(t)\|^{2}_{L_{2}^{\frac{n}{2}}\left(\mathbb{B}\right)}\right)\\[2mm]
&+\frac{2\varepsilon (p-1)}{p}\|g(x)^\frac{1}{p}u(t)\|_{L_{p}^{\frac{n}{p}}\left(\mathbb{B}\right)}^{p}
-2\varepsilon \int_{\mathbb{B}}|u_t(t)|^{m-2}u_t u\frac{dx_{1}}{x_{1}}dx'-2\varepsilon d_1+2\varepsilon H(t).
   \end{split}
\end{equation}
Since we already have $u(t)\in \mathcal{V},\,\,t\geq0$, then by Lemma \ref{k}, we obtain
\begin{equation} \label{5.6}
   \begin{split}\|\nabla_{\mathbb{B}}u\|^{2}_{L_{2}^{\frac{n}{2}}\left(\mathbb{B}\right)}
   -\gamma\|V(x)^{\frac{1}{2}}u\|^{2}_{L_{2}^{\frac{n}{2}}\left(\mathbb{B}\right)}
   &\geq(1-\gamma(C^\ast)^{2})\|\nabla_{\mathbb{B}}u\|^{2}_{L_{2}^{\frac{n}{2}}\left(\mathbb{B}\right)}\\[2mm]
   &>[1-\gamma(C^\ast)^{2}]^\frac{p}{p-2}C_\ast^\frac{-2p}{p-2}\\[2mm]
   &=\frac{2p}{p-2}d,
   \end{split}
\end{equation}
which implies that
$$\frac{2p}{p-2}d<\|\nabla_{\mathbb{B}}u\|^{2}_{L_{2}^{\frac{n}{2}}\left(\mathbb{B}\right)}
-\gamma\|V(x)^{\frac{1}{2}}u\|^{2}_{L_{2}^{\frac{n}{2}}\left(\mathbb{B}\right)}
<\|g(x)^\frac{1}{p}u\|_{L_{p}^{\frac{n}{p}}\left(\mathbb{B}\right)}^{p}.$$
Thus, \begin{equation*}
   \begin{split}&2d_1+(\|\nabla_{\mathbb{B}}u\|^{2}_{L_{2}^{\frac{n}{2}}\left(\mathbb{B}\right)}
   -\gamma\|V(x)^{\frac{1}{2}}u\|^{2}_{L_{2}^{\frac{n}{2}}\left(\mathbb{B}\right)})\\[2mm]
   < &\frac{p-2}{p}\frac{d_1}{d}\|g(x)^\frac{1}{p}u\|_{L_{p}^{\frac{n}{p}}\left(\mathbb{B}\right)}^{p}
   + \|g(x)^\frac{1}{p}u\|_{L_{p}^{\frac{n}{p}}\left(\mathbb{B}\right)}^{p}   \\[2mm]
   =& (1+\frac{p-2}{p}\frac{d_1}{d})\|g(x)^\frac{1}{p}u\|_{L_{p}^{\frac{n}{p}}\left(\mathbb{B}\right)}^{p}.
   \end{split}
\end{equation*}
Combining \eqref{5.5} and \eqref{5.6}, we have
\begin{equation} \label{5.7}
   \begin{split}L'(t)&\geq(1-\eta)H^{-\eta}(t)H'(t)+3\varepsilon\|u_t\|^2_{L^{\frac{n}{2}}_{2}\left(\mathbb{B}\right)}+2\varepsilon H(t)\\[2mm]
&+\varepsilon \frac{p-2}{p}(1-\frac{d_1}{d}) \|g(x)^\frac{1}{p}u\|_{L_{p}^{\frac{n}{p}}\left(\mathbb{B}\right)}^{p}-2\varepsilon \int_{\mathbb{B}}|u_t|^{m-2}u_t u\frac{dx_{1}}{x_{1}}dx'.
   \end{split}
\end{equation}
Hence, from \eqref{5.7}, we can repeat the process after formula \eqref{5.2} in the proof of Theorem \ref{5.1}, we can see the solution of  problem \eqref{1.1} blows up at a finite time.

\textbf{Case 2: $E(0)=d$.}\\
It is  clear from $(u_0, u_1)>0$ that $\|u_1\|^2_{L^{\frac{n}{2}}_{2}\left(\mathbb{B}\right)}\neq 0$. Since $I(u_0)<0$ and
\begin{equation*}
  d=E(0)=\frac12\|u_1\|^2_{L^{\frac{n}{2}}_{2}\left(\mathbb{B}\right)}+J(u_0),
\end{equation*}
we have $u_0\in \mathcal{V}$. Then, we have $u(t)\in \mathcal{V}$, $t\geq0$ by Lemma \ref{j}-(ii). By the similar argument as the proof of Lemma \ref{j}-(ii), we can also have
$\int_{0}^{t_1}\|u_{\tau}\|^{m}_{L^{\frac{n}{m}}_{m}\left(\mathbb{B}\right)}d\tau > 0$ for some small $t_1\in (0,T_{max})$. Then,
it follows that
$$E(t_1)=E(0)-\int_{0}^{t_1}\|u_{\tau}\|^{m}_{L^{\frac{n}{m}}_{m}\left(\mathbb{B}\right)}d\tau
=d-\int_{0}^{t_1}\|u_{\tau}\|^{m}_{L^{\frac{n}{m}}_{m}\left(\mathbb{B}\right)}d\tau < d.$$
Taking $t_1$ as the initial time, it follows the Case $1$ that $u(t)$ blows up at a  finite time. The proof is complete.
\end{proof}

\begin{theorem} \label{z}(Blow up for $E(0)\geq 0$)\\
Let $\gamma\in[0,\frac{1}{(C^{\ast})^{2}})$, $u_{0}\in H_{2,0}^{1,\frac{n}{2}}(\mathbb{B})$,  $u_1\in L^{\frac{n}{2}}_{2}\left(\mathbb{B}\right)$ and $u(t)$ be the  local solution to problem \eqref{1.1} on $[0,T_{max})$. If the initial data satisfies
$$\int_{\mathbb{B}}u_0u_1\frac{dx_{1}}{x_{1}}dx'>\frac{(m-1)M^\frac{1}{m-1}}{m}E(0)\geq0,$$ then the solution  blows up at a finite time, where $M$ is the root of equation $$\frac{K(M)}{\eta(M)}=\frac{(m-1)M^\frac{1}{m-1}}{m}$$ on the interval
$ \left( \frac{(p-m)\alpha+(m-2)(1-\gamma(C^\ast)^{2})\lambda_1}{\lambda_1 (p-2)^2 \alpha (1-\gamma(C^\ast)^{2}) }, + \infty \right) $, and
$$K \left( M \right) = p  - \frac { m - 2 } {\alpha ( p - 2 ) M }, \eta  ( M ) = \sqrt { ( 2 + K ( M ) ) \left[ ( K ( M ) - 2 ) \lambda_{ 1 }(1-\gamma(C^\ast)^{2}) - \frac { p - m } { ( p - 2 ) M } \right] },$$
where $\alpha=\inf_{x\in \mathbb{B}}g(x)>0$ and  $\lambda_1$ is the first eigenvalue of the operator $-\Delta_{\mathbb{B}}$  given in Lemma \ref{c}.
\end{theorem}
\begin{proof}
Arguing by contradiction, we assume that the maximal existence time $T_{max}=\infty$. Suppose $u(t)$ is the global solution, for $t\in[0,+\infty)$, setting
$$ L_{1}(t)=F^{\prime}(t)-\frac{2(m-1)M^{\frac{1}{m-1}}}{m}E(t),$$
where $F(t)=\|u\|^2_{L^{\frac{n}{2}}_{2}\left(\mathbb{B}\right)}$.
In view of  Young's inequality, we have
$$\int_{\mathbb{B}}|u_t|^{m-2}u_t u\frac{dx_{1}}{x_{1}}dx'
\leq \frac{m-1}{m}\varepsilon^{-\frac{m}{m-1}}\|u_t\|^m_{L^{\frac{n}{m}}_{m}\left(\mathbb{B}\right)}
+\frac{\varepsilon^{m}}{m}\|u\|^m_{L^{\frac{n}{m}}_{m}\left(\mathbb{B}\right)}$$ for any $\varepsilon>0$.
Setting $\varepsilon^{-\frac{m}{m-1}}=M^{\frac{1}{m-1}}$, we get
\begin{equation*}
   \begin{split}L_1'(t)&=2\|u_t\|^2_{L^{\frac{n}{2}}_{2}\left(\mathbb{B}\right)}
   -2\|\nabla_{\mathbb{B}}u\|^{2}_{L_{2}^{\frac{n}{2}}\left(\mathbb{B}\right)}
   +2\gamma\|V(x)^{\frac{1}{2}}u\|^{2}_{L_{2}^{\frac{n}{2}}\left(\mathbb{B}\right)}
   +2\|g(x)^\frac{1}{p}u\|_{L_{p}^{\frac{n}{p}}\left(\mathbb{B}\right)}^{p}\\[2mm]
   &-2\int_{\mathbb{B}}|u_t|^{m-2}u_t u\frac{dx_{1}}{x_{1}}dx'
   +2\frac{(m-1)}{m}M^\frac{1}{m-1}\|u_t\|^{m}_{L^{\frac{n}{m}}_{m}\left(\mathbb{B}\right)}\\[2mm]
   &\geq 2\|u_t\|^2_{L^{\frac{n}{2}}_{2}\left(\mathbb{B}\right)}
   -2\|\nabla_{\mathbb{B}}u\|^{2}_{L_{2}^{\frac{n}{2}}\left(\mathbb{B}\right)}
   +2\gamma\|V(x)^{\frac{1}{2}}u\|^{2}_{L_{2}^{\frac{n}{2}}\left(\mathbb{B}\right)}
   +2\|g(x)^\frac{1}{p}u\|_{L_{p}^{\frac{n}{p}}\left(\mathbb{B}\right)}^{p}\\[2mm]
   &-\frac{2}{Mm}\|u\|^{m}_{L^{\frac{n}{m}}_{m}\left(\mathbb{B}\right)}
   -2\frac{(m-1)}{m}M^\frac{1}{m-1}\|u_t\|^{m}_{L^{\frac{n}{m}}_{m}\left(\mathbb{B}\right)}
   +2\frac{(m-1)}{m}M^\frac{1}{m-1}\|u_t\|^{m}_{L^{\frac{n}{m}}_{m}\left(\mathbb{B}\right)}\\[2mm]
   &=2\|u_t\|^2_{L^{\frac{n}{2}}_{2}\left(\mathbb{B}\right)}
   -2\|\nabla_{\mathbb{B}}u\|^{2}_{L_{2}^{\frac{n}{2}}\left(\mathbb{B}\right)}
   +2\gamma\|V(x)^{\frac{1}{2}}u\|^{2}_{L_{2}^{\frac{n}{2}}\left(\mathbb{B}\right)}
   +2\|g(x)^\frac{1}{p}u\|_{L_{p}^{\frac{n}{p}}\left(\mathbb{B}\right)}^{p}
   -\frac{2}{Mm}\|u\|^{m}_{L^{\frac{n}{m}}_{m}\left(\mathbb{B}\right)}.
   \end{split}
\end{equation*}

Since the function $\frac{x^y}{y}(x>0,y>0)$ is convex with respect to $y$, we arrive at
$$\frac{1}{m}\|u\|^{m}_{L^{\frac{n}{m}}_{m}\left(\mathbb{B}\right)}
\leq\frac{m-2}{(p-2)p}\|u\|^{p}_{L^{\frac{n}{p}}_{p}\left(\mathbb{B}\right)}+
\frac{p-m}{2(p-2)}\|u\|^{2}_{L^{\frac{n}{2}}_{2}\left(\mathbb{B}\right)}.$$
Noticing $0<\alpha=\displaystyle\inf_{x\in\mathbb{B}}g(x)$, we have
\begin{equation*}
   \begin{split}L_1'(t)&\geq2\|u_t\|^2_{L^{\frac{n}{2}}_{2}\left(\mathbb{B}\right)}
   -2\|\nabla_{\mathbb{B}}u\|^{2}_{L_{2}^{\frac{n}{2}}\left(\mathbb{B}\right)}
   +2\gamma\|V(x)^{\frac{1}{2}}u\|^{2}_{L_{2}^{\frac{n}{2}}\left(\mathbb{B}\right)}
   +2\|g(x)^\frac{1}{p}u\|_{L_{p}^{\frac{n}{p}}\left(\mathbb{B}\right)}^{p}\\[2mm]
   &-\frac{2}{M}\frac{m-2}{p(p-2)}\|u\|_{L_{p}^{\frac{n}{p}}\left(\mathbb{B}\right)}^{p}
   -\frac{p-m}{M(p-2)}\|u\|_{L_{2}^{\frac{n}{2}}\left(\mathbb{B}\right)}^{2}\\[2mm]
   &\geq2\|u_t\|^2_{L^{\frac{n}{2}}_{2}\left(\mathbb{B}\right)}
   -2(\|\nabla_{\mathbb{B}}u\|^{2}_{L_{2}^{\frac{n}{2}}\left(\mathbb{B}\right)}
   -\gamma\|V(x)^{\frac{1}{2}}u\|^{2}_{L_{2}^{\frac{n}{2}}\left(\mathbb{B}\right)})
   \\[2mm]
   &+\frac{2}{p}(p-\frac{m-2}{M\alpha(p-2)})\|g(x)^\frac{1}{p}u\|_{L_{p}^{\frac{n}{p}}\left(\mathbb{B}\right)}^{p}
   -\frac{p-m}{M(p-2)}\|u\|_{L_{2}^{\frac{n}{2}}\left(\mathbb{B}\right)}^{2}.
   \end{split}
\end{equation*}
By Lemma \ref{c} and  the definition of $E(t)$, we have
\begin{equation*}
   \begin{split}L_1'(t)\geq&(2+K(M))\|u_t\|^2_{L^{\frac{n}{2}}_{2}\left(\mathbb{B}\right)}
   -\frac{p-m}{M(p-2)}\|u\|_{L_{2}^{\frac{n}{2}}\left(\mathbb{B}\right)}^{2}-2K(M)E(t)\\[2mm]
   &+(K(M)-2)(\|\nabla_{\mathbb{B}}u\|^{2}_{L_{2}^{\frac{n}{2}}\left(\mathbb{B}\right)}
   -\gamma\|V(x)^{\frac{1}{2}}u\|^{2}_{L_{2}^{\frac{n}{2}}\left(\mathbb{B}\right)})
   \\[2mm]
   \geq &(2+K(M))\|u_t\|^2_{L^{\frac{n}{2}}_{2}\left(\mathbb{B}\right)}
   +(K(M)-2)(1-\gamma(C^\ast)^{2})\|\nabla_{\mathbb{B}}u\|^{2}_{L_{2}^{\frac{n}{2}}\left(\mathbb{B}\right)}\\[2mm]
   &-\frac{p-m}{M(p-2)}\|u\|_{L_{2}^{\frac{n}{2}}\left(\mathbb{B}\right)}^{2}-2K(M)E(t)\\[2mm]
   \geq& (2+K(M))\|u_t\|^2_{L^{\frac{n}{2}}_{2}\left(\mathbb{B}\right)}
   +(K(M)-2)(1-\gamma(C^\ast)^{2})\lambda_1\|u\|^2_{L^{\frac{n}{2}}_{2}\left(\mathbb{B}\right)}\\[2mm]
   &-\frac{p-m}{M(p-2)}\|u\|_{L_{2}^{\frac{n}{2}}\left(\mathbb{B}\right)}^{2}-2K(M)E(t),
   \end{split}
\end{equation*}
where $ K ( M ) = p  - \frac { m - 2 } { \alpha( p - 2 ) M } > 2.$  Denote by
$$ K_{ 1 } ( M ) :=\lambda _ { 1 } ( K ( M ) - 2 )(1-\gamma(C^\ast)^{2}) - \frac { p - m } { ( p - 2 ) M } > 0,$$
which implies $$L_1'(t)\geq(2+K(M))\|u_t\|^2_{L^{\frac{n}{2}}_{2}\left(\mathbb{B}\right)}+K _ { 1 } ( M )\|u\|_{L_{2}^{\frac{n}{2}}\left(\mathbb{B}\right)}^{2}-2K(M)E(t).$$
By Cauchy inequality and cone H\"{o}lder inequality, we get
$$ (2+K(M))\|u_{t}\|^{2}_{L^{\frac{n}{2}}_{2}\left(\mathbb{B}\right)}+K_{1}(M) \|u\|^{2}_{L^{\frac{n}{2}}_{2}\left(\mathbb{B}\right)} \geq 2 \sqrt { (2 + K ( M ))   K _ { 1 } ( M )  }( u , u _ { t } ) .$$
Thus,
\begin{equation}\label{5.8}
  L_ { 1 } ^ { \prime } ( t ) \geq \eta ( M ) ( 2 ( u , u _ { t } ) - \frac { 2 K ( M ) } { \eta ( M ) } E ( t ) ) ,
\end{equation}
where $$ \eta( M ) = \sqrt { ( 2 + K ( M ) ) K _ { 1 } ( M ) } .$$
Choosing  $ M _ { 0 } = \frac { ( m - 2 ) \lambda _ { 1 }(1-\gamma(C^\ast)^{2}) + (p - m) \alpha} { ( p - 2 ) ^ { 2 } \lambda _ { 1 }\alpha (1-\gamma(C^\ast)^{2})},$
we easily obtain
$$\lim\limits_ { M \to  M _ { 0 } } \frac { K ( M ) } { \eta ( M ) } = + \infty ,\lim _ { M \rightarrow M _ { 0 } } \frac {  (m-1) M ^ { \frac { 1 } { m-1 } } } { m } = \frac { (m-1) M _ { 0 } ^ { \frac { 1 } { m-1 } } } { m  } ,$$
$$ \lim _ { M \to + \infty } \frac { K ( M ) } { \eta ( M ) } = \frac { p } { \sqrt { ( p + 2 ) ( p - 2 ) (1-\gamma(C^\ast)^{2})\lambda _ { 1 } } } , \lim _ { M \to + \infty } \frac { (m-1) M ^ { \frac { 1 } {m-1 } } } { m  } = + \infty ,$$
Letting  $$\varphi(M)=\frac { K ( M ) } { \eta ( M ) } - \frac { (m-1) M ^ { \frac { 1 } { m-1 } } } { m },$$ it is clear that $\lim\limits_ { M \to  M _ { 0 } }\varphi(M)=+\infty$ and $\lim\limits_ { M \to  +\infty }\varphi(M)=-\infty$.\\
Thus, there must exist  $M>M_0$,  such that
\begin{equation}\label{5.9}
  \frac { K ( M ) } { \eta ( M ) } = \frac { (m-1) M ^ { \frac { 1 } { m-1 } } } { m } .
\end{equation}
Choosing $M$  to satisfy the \eqref{5.9}, we get
\begin{equation*}
 L_{1}^{\prime}(t)\geq\eta( M ) L _ { 1 } ( t ),\,\, L _ { 1 } ( 0 ) = 2\int_{\mathbb{B}}u_0 u_1 \frac{dx_{1}}{x_{1}}dx' - 2\frac { (m-1) M ^ { \frac{1}{m-1}} } { m  } E ( 0 ) > 0.
\end{equation*}
Hence, we obtain
\begin{equation*}
L _ { 1 } ( t ) \geq L _ { 1 } ( 0 ) e ^ { \eta ( M ) t },\,\,  t \geq 0.
\end{equation*}
Since $u(x,t)$ is the global solution, we can get $0\leq E(t)\leq E(0)$ for all  $t\in [0,+\infty)$ ( otherwise, if there exists a $t'\in (0, \infty)$ such that $E(t')<0$, we can obtain the solution would blow up at a finite time by Theorem \ref{th5.1}. Therefore,
\begin{equation}\label{5.10}
  F ^ { \prime } ( t )  \geq L _ { 1 } ( 0 ) e ^ { \eta  ( M ) t } , \,\,t \geq 0,
\end{equation}
which implies
\begin{equation}\label{5.11}
 F(t)=\|u\|^2_{L^{\frac{n}{2}}_{2}\left(\mathbb{B}\right)}\geq \|u_0\|^2_{L^{\frac{n}{2}}_{2}\left(\mathbb{B}\right)} + \frac { L _ { 1 }  ( 0 )} { \eta ( M ) }  \left( e ^ { \eta ( M ) t } - 1 \right),\,\,\,\,\, t\geq0.
\end{equation}

On the other hand, from H\"{o}lder inequality, we have
\begin{equation} \label{5.12}
   \begin{split}\|u\|_{L^{\frac{n}{2}}_{2}\left(\mathbb{B}\right)}&\leq\|u_0\|_{L^{\frac{n}{2}}_{2}\left(\mathbb{B}\right)}+
   \int_{0}^{t}\|u_{\tau}\|_{L^{\frac{n}{2}}_{2}\left(\mathbb{B}\right)}d\tau\\[2mm]
   &\leq\|u_0\|_{L^{\frac{n}{2}}_{2}\left(\mathbb{B}\right)}+
   C(\int_{0}^{t}\|u_{\tau}\|^m_{L^{\frac{n}{m}}_{m}\left(\mathbb{B}\right)}d\tau)^\frac{1}{m}t^{\frac{m-1}{m}}\\[2mm]
   &=\|u_0\|_{L^{\frac{n}{2}}_{2}\left(\mathbb{B}\right)}+
   Ct^{\frac{m-1}{m}}[E(0)-E(t)]^\frac{1}{m}\\[2mm]
   &\leq \|u_0\|_{L^{\frac{n}{2}}_{2}\left(\mathbb{B}\right)}+
   Ct^{\frac{m-1}{m}}[E(0)]^\frac{1}{m}.
   \end{split}
\end{equation}
Then the estimate \eqref{5.12} contradicts the inequality \eqref{5.11}.
Thus, the solution of  problem \eqref{1.1}  cannot be extended to the whole interval $[0, \infty)$. The proof is complete.
\end{proof}

In view of  Theorem \ref{z}, we can prove that the existence of certain  solutions, which blow up with arbitrary initial energy
level (including the case  $E(0) > d$).

\begin{corollary}
For any given constant $R$ including $R>d$, there exists two function $u^R_0\in \mathcal{H}_{2,0}^{1,\frac{n}{2}}
\left(\mathbb{B}\right)$
and $u^R_1\in L^\frac{n}{2}_2(\mathbb{B})$ such that $E(0)=R$, and the corresponding  solution $u(t)$ to  problem \eqref{1.1} with initial data $u_0= u^R_0$, $u_1= u^R_1 $  blows up at a finite time, where
$$E(0)=\frac{1}{2}\|u^R_1\|^2_{L^{\frac{n}{2}}_{2}\left(\mathbb{B}\right)}
+\frac{1}{2}\|\nabla_{\mathbb{B}}u^R_0\|^{2}_{L_{2}^{\frac{n}{2}}\left(\mathbb{B}\right)}
-\frac{1}{2}\gamma\|V(x)^{\frac{1}{2}}u^R_0\|^{2}_{L_{2}^{\frac{n}{2}}\left(\mathbb{B}\right)}
-\frac{1}{p}\|g(x)^{\frac{1}{p}}u^R_0\|^{p}_{L^{\frac{n}{p}}_{p}\left(\mathbb{B}\right)}.$$
\end{corollary}
\begin{proof}
For any nonzero functions $\omega_1(x)\in \mathcal{H}_{2,0}^{1,\frac{n}{2}}
\left(\mathbb{B}\right)$, $\omega_2(x)\in L^\frac{n}{2}_2(\mathbb{B}) $ satisfying
$(\omega_1,\omega_2)=0$, we choose the initial data $u^R_0$ and $u^R_1$ as $$u^R_0=r_1\omega_1(x),\,\,u^R_1=r_1\omega_1(x)+r_2\omega_2(x),$$ where $r_1$ and $r_2$ are two positive constants to be chosen  later. Let $u(t)$ be the corresponding  solution to  problem \eqref{1.1} with initial data
$u_0= u^R_0$ and $u_1= u^R_1 $. Then, the initial energy can be writen as
\begin{equation*}
   \begin{split}E(0)   &= \frac{1}{2}r_1^2\|\omega_1(x)\|^2_{L^{\frac{n}{2}}_{2}\left(\mathbb{B}\right)}
   +\frac{1}{2}r_2^2\|\omega_2(x)\|^2_{L^{\frac{n}{2}}_{2}\left(\mathbb{B}\right)}
   +\frac{1}{2}r_1^2\|\nabla_{\mathbb{B}}\omega_1(x)\|^{2}_{L_{2}^{\frac{n}{2}}\left(\mathbb{B}\right)}\\[2mm]
   &-\frac{1}{2}r_1^2\gamma\|V(x)^{\frac{1}{2}}\omega_1(x)\|^{2}_{L_{2}^{\frac{n}{2}}\left(\mathbb{B}\right)}
   -\frac{1}{p}r_1^p\|g(x)^{\frac{1}{p}}\omega_1(x)\|^{p}_{L^{\frac{n}{p}}_{p}\left(\mathbb{B}\right)}\\[2mm]
   &=\frac{1}{2}r_2^2\|\omega_2(x)\|^2_{L^{\frac{n}{2}}_{2}\left(\mathbb{B}\right)}+\chi(r_1),
   \end{split}
\end{equation*}
where$$\chi(r_1)=\frac{1}{2}r_1^2\|\omega_1(x)\|^2_{L^{\frac{n}{2}}_{2}\left(\mathbb{B}\right)}
   +\frac{1}{2}r_1^2\|\nabla_{\mathbb{B}}\omega_1(x)\|^{2}_{L_{2}^{\frac{n}{2}}\left(\mathbb{B}\right)}
-\frac{1}{2}r_1^2\gamma\|V(x)^{\frac{1}{2}}\omega_1(x)\|^{2}_{L_{2}^{\frac{n}{2}}\left(\mathbb{B}\right)}
-\frac{1}{p}r_1^p\|g(x)^{\frac{1}{p}}\omega_1(x)\|^{p}_{L^{\frac{n}{p}}_{p}\left(\mathbb{B}\right)}.$$
Since $p>2$, we can see that
$\lim\limits_{r_1\to+\infty}\chi(r_1)=-\infty$. Therefore, we can take  $r_1$ sufficiently
large such that $$\chi(r_1)<R<\frac{m}{(m-1)M^\frac{1}{m-1}}(u^R_0,u^R_1)=r_1^2\frac{m}{(m-1)M^\frac{1}{m-1}}\|\omega_1(x)\|^{2}_{L_{2}^{\frac{n}{2}}\left(\mathbb{B}\right)}.$$
For the selected $r_1$, taking  $r_2=\frac{\sqrt{2(R-\chi(r_1))}}{\|\omega_2(x)\|^{2}_{L_{2}^{\frac{n}{2}}\left(\mathbb{B}\right)}}$, we have
$$E(0)=\frac{2(R-\chi(r_1))}{\|\omega_2(x)\|^{2}_{L_{2}^{\frac{n}{2}}\left(\mathbb{B}\right)}}\|\omega_2(x)\|^{2}_{L_{2}^{\frac{n}{2}}\left(\mathbb{B}\right)}+\chi(r_1)=R.$$
Thus,  it follows Theorem \ref{z} that the corresponding solution $u(t)$ to  problem \eqref{1.1} with initial data $u^R_0\in \mathcal{H}_{2,0}^{1,\frac{n}{2}}
\left(\mathbb{B}\right)$
and $u^R_1\in L^\frac{n}{2}_2(\mathbb{B})$  blows up at a finite time.  The proof is complete.
\end{proof}

From Theorem \ref{th5.1} and \ref{th5.2}, we have the lifespan or the upper bound of the blowup time provided the initial energy are
at both subcritical and critical levels. We are in the position to derive an upper bound of the blowup time at high energy level with some additional assumption.

\begin{theorem} (an upper bound of the blow-up time)\\
Let all assumptions in Theorem $5.3$ be fulfilled,  and
$$\|u_0\|^2_{L^{\frac{n}{2}}_{2}\left(\mathbb{B}\right)}\geq\frac{p+2+\xi}{2\rho_1}E(0),$$
where $\xi$ is any positive number and $\rho_1$ is shown as \eqref{5.17}.  Then the solution $u(t)$ to problem \eqref{1.1} blows up at some finite time $T_{max}$ in the sense that  $\lim\limits_{t\to T_{max}}\|u\|_{L_{p}^{\frac{n}{p}}\left(\mathbb{B}\right)}= +\infty$
and the blow up time $T_{max}$ can be estimated from above as follows
$$ T _{max} \leq \mathcal{F} ^ { - \frac { \sigma } { 1 - \sigma } } ( 0 ) \frac {\mathcal{M} _ { 1 } } { \mathcal{M} _ { 2 } } \frac { 1 - \sigma } { \sigma } ,$$
where $ 0 < \sigma \leq \frac { p  - 2 } { 2 p  } $ and
$\mathcal{M}_1$ , $\mathcal{M}_2$ will be  presented in \eqref{5.19} and \eqref{5.21} respectively.
$\mathcal{F}(0)$ can be found in \eqref{5.22}.
\end{theorem}
\begin{proof}
Without loss of generality, we may assume $E(t)\geq 0 $ for all $t\geq0$. Define
\begin{equation}\label{5.13}
 \mathcal{H}(t)=E(0)-E(t)
\end{equation}
for $t\geq 0.$
Let
$$\mathcal{F} ( t )=\mathcal{H}^{1-\sigma} ( t )+\varepsilon\int_{\mathbb{B}}u_t u\frac{dx_{1}}{x_{1}}dx',$$
where $\varepsilon > 0$ will be determined  later. The proof will be divided into three
steps.\\
\textbf{Step 1: Estimate for $\mathcal{F} '( t )$.}\\
 Differentiating  directly  $\mathcal{F}( t )$,  for $0 <k< 1$, we have
\begin{equation} \label{5.14}
   \begin{split}\mathcal{F} '( t )&=(1-\sigma){\mathcal{H}}^{-\sigma} ( t ) \mathcal{H}'(t)
   +\varepsilon \|u_t\|^2_{L^{\frac{n}{2}}_{2}\left(\mathbb{B}\right)}
   -\varepsilon\|\nabla_{\mathbb{B}}u\|^{2}_{L_{2}^{\frac{n}{2}}\left(\mathbb{B}\right)} \\[2mm] &+\varepsilon\gamma\|V(x)^{\frac{1}{2}}u\|^{2}_{L_{2}^{\frac{n}{2}}\left(\mathbb{B}\right)}
   -\varepsilon\int_{\mathbb{B}}|u_t|^{m-2}u_t u\frac{dx_{1}}{x_{1}}dx'+k\varepsilon\|g(x)^\frac{1}{p}u\|_{L_{p}^{\frac{n}{p}}\left(\mathbb{B}\right)}^{p}\\[2mm]
   &+\varepsilon p (1-k)(\frac{1}{2}\|u_t\|^2_{L^{\frac{n}{2}}_{2}\left(\mathbb{B}\right)}
   +\frac{1}{2}\|\nabla_{\mathbb{B}}u\|^{2}_{L_{2}^{\frac{n}{2}}\left(\mathbb{B}\right)}
   -\frac{1}{2}\gamma\|V(x)^{\frac{1}{2}}u\|^{2}_{L_{2}^{\frac{n}{2}}\left(\mathbb{B}\right)}-E(t)).   \end{split}
\end{equation}
By Young's inequality with $\varepsilon_2 > 1$, $\mathcal{H}'(t)=-E'(t)$, we have
\begin{equation} \label{5.15}
   \begin{split}&|\int_{\mathbb{B}}|u_t|^{m-2}u_t u\frac{dx_{1}}{x_{1}}dx'|\\[2mm]
   &\leq\int_{\mathbb{B}}|u_t|^{m-1}\mathcal{H}^{-\frac{\sigma(m-1)}{m}}u \mathcal{H}^{\frac{\sigma(m-1)}{m}}\frac{dx_{1}}{x_{1}}dx'\\[2mm]
   &\leq \varepsilon_2 \frac{m-1}{m}\int_{\mathbb{B}}{\mathcal{H}}^{-\sigma} ( t )|u_t|^{m}\frac{dx_{1}}{x_{1}}dx'
   +\frac{1}{m}\frac{1}{\varepsilon_2^{m-1}}\int_{\mathbb{B}}{\mathcal{H}}^{\sigma(m-1)} ( t )|u|^{m}\frac{dx_{1}}{x_{1}}dx'\\[2mm]
   &\leq \varepsilon_2{\mathcal{H}}^{-\sigma} ( t ){\mathcal{H}}'(t)
   +\frac{E^{\sigma(m-1)}(0)}{\varepsilon_2^{m-1}}
   (\frac{p-m}{p-2}\|u\|^2_{L^{\frac{n}{2}}_{2}\left(\mathbb{B}\right)}
   +\frac{m-2}{p-2}\|u\|_{L_{p}^{\frac{n}{p}}\left(\mathbb{B}\right)}^{p}).
   \end{split}
\end{equation}
Inserting \eqref{5.15} into \eqref{5.14}, we  can see that
\begin{equation} \label{5.16}
   \begin{split}\mathcal{F}'(t)\geq&(1-\sigma){\mathcal{H}}^{-\sigma} ( t ) \mathcal{H}'(t)
   +\varepsilon \|u_t\|^2_{L^{\frac{n}{2}}_{2}\left(\mathbb{B}\right)}
   -\varepsilon(\|\nabla_{\mathbb{B}}u\|^{2}_{L_{2}^{\frac{n}{2}}\left(\mathbb{B}\right)}
   -\gamma\|V(x)^{\frac{1}{2}}u\|^{2}_{L_{2}^{\frac{n}{2}}\left(\mathbb{B}\right)})\\[2mm]
   &- \varepsilon\varepsilon_2{\mathcal{H}}^{-\sigma} ( t ) \mathcal{H}'(t)
   -\varepsilon\frac{1}{\varepsilon_2^{m-1}}E^{\sigma(m-1)}(0)\frac{p-m}{p-2}\|u\|^2_{L^{\frac{n}{2}}_{2}\left(\mathbb{B}\right)}
   -\varepsilon\frac{m-2}{p-2}E^{\sigma(m-1)}(0)\|u\|_{L_{p}^{\frac{n}{p}}\left(\mathbb{B}\right)}^{p}\\[2mm]
   &+\varepsilon p (1-k)[\frac{1}{2}\|u_t\|^2_{L^{\frac{n}{2}}_{2}\left(\mathbb{B}\right)}
   +\frac{1}{2}(\|\nabla_{\mathbb{B}}u\|^{2}_{L_{2}^{\frac{n}{2}}\left(\mathbb{B}\right)}
   -\gamma \|V(x)^{\frac{1}{2}}u\|^{2}_{L_{2}^{\frac{n}{2}}\left(\mathbb{B}\right)})-E(t)]\\[2mm]
   \geq& (1-\sigma-\varepsilon\varepsilon_2){\mathcal{H}}^{-\sigma} ( t ) \mathcal{H}'(t)
   +\varepsilon(\frac{p(1-k)+2}{2})\|u_t\|^2_{L^{\frac{n}{2}}_{2}\left(\mathbb{B}\right)}\\[2mm]
   &+\varepsilon[\frac{p(1-k)-2}{2}](\|\nabla_{\mathbb{B}}u\|^{2}_{L_{2}^{\frac{n}{2}}\left(\mathbb{B}\right)}
   -\gamma \|V(x)^{\frac{1}{2}}u\|^{2}_{L_{2}^{\frac{n}{2}}\left(\mathbb{B}\right)})-\varepsilon p(1-k)E(t)\\[2mm]
   &-\varepsilon\frac{1}{\varepsilon_2^{m-1}}E^{\sigma(m-1)}(0)\frac{p-m}{p-2}\|u\|^2_{L^{\frac{n}{2}}_{2}\left(\mathbb{B}\right)}
   +\varepsilon (k\alpha-\frac{1}{\varepsilon_2^{m-1}}E^{\sigma(m-1)}(0))\|u\|_{L_{p}^{\frac{n}{p}}\left(\mathbb{B}\right)}^{p}.
   \end{split}
\end{equation}
Let us choose $k=\frac{p-2}{2p}>0$, use the embedding theorem and recall \eqref{5.13}, then \eqref{5.16} can be transferred to
\begin{equation*}
   \begin{split}\mathcal{F}'(t)&\geq (1-\sigma-\varepsilon\varepsilon_2){\mathcal{H}}^{-\sigma} ( t ) \mathcal{H}'(t)
   +\varepsilon(\frac{p+6}{4})\|u_t\|^2_{L^{\frac{n}{2}}_{2}\left(\mathbb{B}\right)}
   -\varepsilon(\frac{p+2}{2})E(0)+\varepsilon(\frac{p+2}{2})\mathcal{H}(t)\\[2mm]
   &-\varepsilon E^{\sigma(m-1)}(0) \frac{p-m}{(p-2)\varepsilon_2^{m-1}}\|u\|^2_{L^{\frac{n}{2}}_{2}\left(\mathbb{B}\right)}
   +\varepsilon[\frac{\alpha(p-2)}{2p}-\frac{1}{\varepsilon_2^{m-1}}E^{\sigma(m-1)}(0)]\|u\|_{L_{p}^{\frac{n}{p}}\left(\mathbb{B}\right)}^{p}\\[2mm]
   &+\varepsilon\frac{p-2}{4}(\|\nabla_{\mathbb{B}}u\|^{2}_{L_{2}^{\frac{n}{2}}\left(\mathbb{B}\right)}-
   \gamma\|V(x)^{\frac{1}{2}}u\|^{2}_{L_{2}^{\frac{n}{2}}\left(\mathbb{B}\right)})\\[2mm]
   &\geq (1-\sigma-\varepsilon\varepsilon_2){\mathcal{H}}^{-\sigma} ( t ) \mathcal{H}'(t)
   +\varepsilon(\frac{p+6}{4})\|u_t\|^2_{L^{\frac{n}{2}}_{2}\left(\mathbb{B}\right)}
   -\varepsilon(\frac{p+2}{2})E(0)\\[2mm]
   &+\varepsilon(\frac{p+2}{2})\mathcal{H}(t)
   +\varepsilon[\frac{\lambda_1^2(1-\gamma(C^\ast)^{2})(p-2)}{4}- E^{\sigma(m-1)}(0) \frac{p-m}{(p-2)\varepsilon_2^{m-1}}]\|u\|^2_{L^{\frac{n}{2}}_{2}\left(\mathbb{B}\right)}\\[2mm]
   &+\varepsilon(\frac{\alpha(p-2)}{2p}-\frac{1}{\varepsilon_2^{m-1}}E^{\sigma(m-1)}(0))\|u\|_{L_{p}^{\frac{n}{p}}\left(\mathbb{B}\right)}^{p}.
   \end{split}
\end{equation*}
Next, we can choose $\varepsilon_2 > 1$ such that
\begin{equation}\label{5.17}
\rho_1:=\frac{\lambda_1^2(1-\gamma(C^\ast)^{2})(p-2)}{4}- E^{\sigma(m-1)}(0)
\frac{p-m}{(p-2)\varepsilon_2^{m-1}}>0,
\end{equation}
$$\rho_2:=\frac{\alpha(p-2)}{2p}-\frac{1}{\varepsilon_2^{m-1}}E^{\sigma(m-1)}(0)>0.$$
From \eqref{5.10}, we have
\begin{equation*}
   \begin{split}&\int_{\mathbb{B}}u u_t\frac{dx_{1}}{x_{1}}dx'-\frac{m-1}{m}M^{\frac{1}{m-1}}E(t)\\[2mm]
&\geq(\int_{\mathbb{B}}u_0 u_1\frac{dx_{1}}{x_{1}}dx'-\frac{m-1}{m}M^{\frac{1}{m-1}}E(0))e^{\eta(M)t}\\[2mm]
&>0,
\end{split}
\end{equation*}
which implies that
\begin{equation*}
   \begin{split}\frac{d}{dt}\|u\|^2_{L^{\frac{n}{2}}_{2}\left(\mathbb{B}\right)}&=2\int_{\mathbb{B}}u u_t\frac{dx_{1}}{x_{1}}dx'  \\[2mm]
&\geq2[\frac{m-1}{m}M^{\frac{1}{m-1}}E(t)-\int_{\mathbb{B}}u_0 u_1\frac{dx_{1}}{x_{1}}dx'
-\frac{m-1}{m}M^{\frac{1}{m-1}}E(0)]e^{\eta(M)t}\\[2mm]
&\geq0.
\end{split}
\end{equation*}
Thus, we have $$\|u\|^2_{L^{\frac{n}{2}}_{2}\left(\mathbb{B}\right)}\geq\|u_0\|^2_{L^{\frac{n}{2}}_{2}\left(\mathbb{B}\right)}
\geq\frac{p+2+\xi}{2\rho_1}E(0).$$
By choosing $\varepsilon$ sufficiently small such that $1-\sigma-\varepsilon\varepsilon_2 >0 $, it follows that
\begin{equation} \label{5.18}
   \begin{split}\mathcal{F}'(t)&\geq\varepsilon(\frac{p+6}{4})\|u_t\|^2_{L^{\frac{n}{2}}_{2}\left(\mathbb{B}\right)}
   +\varepsilon\rho_1\frac{p+2+\xi}{2\rho_1}E(0)-\varepsilon\frac{p+2}{2}E(0)+\varepsilon\frac{p+2}{2}\mathcal{H}(t)+
   \varepsilon\rho_2\|u\|_{L_{p}^{\frac{n}{p}}\left(\mathbb{B}\right)}^{p}\\[2mm]
   &\geq\mathcal{M}_{1} [\|u\|_{L_{p}^{\frac{n}{p}}\left(\mathbb{B}\right)}^{p}+1+\|u_t\|^2_{L^{\frac{n}{2}}_{2}\left(\mathbb{B}\right)}+\mathcal{H}(t)],
   \end{split}
\end{equation}
where
\begin{equation}\label{5.19}
  \mathcal{M}_1:=\varepsilon\min\left\{\frac{p+6}{4},\frac{p+2}{2},\frac{\xi}{2}E(0),\rho_2\right\}
\end{equation}\\
\textbf{Step 2: Estimate for $\mathcal{F} ^{\frac{1}{1-\sigma}}( t )$.}\\
In what follows, let us consider
$$\mathcal{F}^ {\frac{1}{1-\sigma}}( t )=[\mathcal{H}^{1-\sigma} ( t )+\varepsilon\int_{\mathbb{B}}|u_t|^{m-2}u_t u\frac{dx_{1}}{x_{1}}dx']^{\frac{1}{1-\sigma}}.$$
From the Young's inequality and ${L^{\frac{n}{p}}_{p}\left(\mathbb{B}\right)}\hookrightarrow{L^{\frac{n}{2}}_{2}\left(\mathbb{B}\right)}$ , we have
\begin{equation*}
   \begin{split}{|\int_{\mathbb{B}}u_t u\frac{dx_{1}}{x_{1}}dx'|}^{\frac{1}{1-\sigma}}
   &\leq C \|u_t\|_{L^{\frac{n}{2}}_{2}\left(\mathbb{B}\right)}^{\frac{1}{1-\sigma}} \|u\|_{L_{p}^{\frac{n}{p}}\left(\mathbb{B}\right)}^{\frac{1}{1-\sigma}}\\[2mm]
   &\leq C_1\|u_t\|^2_{L^{\frac{n}{2}}_{2}\left(\mathbb{B}\right)}+C_2\|u\|^{\frac{2}{2(1-\sigma)-1}}_{L_{p}^{\frac{n}{p}}\left(\mathbb{B}\right)}.
   \end{split}
\end{equation*}
Using $0<\sigma<\frac{p-2}{2p}$, we arrive at
$$|\int_{\mathbb{B}}u_t u\frac{dx_{1}}{x_{1}}dx'|\leq C_1\|u_t\|^2_{L^{\frac{n}{2}}_{2}\left(\mathbb{B}\right)}
+C_2 {\frac{2}{p[2(1-\sigma)-1]}}\|u\|_{L_{p}^{\frac{n}{p}}\left(\mathbb{B}\right)}^{p}
+C_2 {\frac{p[2(1-\sigma)-1]-2}{p[2(1-\sigma)-1]}}. $$
On the other hand,
\begin{equation} \label{5.20}
   \begin{split}\mathcal{F}^{\frac{1}{1-\sigma}}(t)&\leq
   2^{{\frac{1}{1-\sigma}}-1}\mathcal{H}( t )+[\varepsilon\int_{\mathbb{B}}|u_t|^{m-2}u_t u\frac{dx_{1}}{x_{1}}dx']^{\frac{1}{1-\sigma}}\\[2mm]
   &\leq  2^{\frac{\sigma}{1-\sigma}}
   [\mathcal{H}(t)+\varepsilon^{\frac{1}{1-\sigma}}C_1\|u_t\|^2_{L^{\frac{n}{2}}_{2}}
   +\varepsilon^{\frac{1}{1-\sigma}}C_2 {\frac{2}{p[2(1-\sigma)-1]}}\|u\|_{L_{p}^{\frac{n}{p}}}^{p}
   +\varepsilon^{\frac{1}{1-\sigma}}C_2 {\frac{p[2(1-\sigma)-1]-2}{p[2(1-\sigma)-1]}}]\\[2mm]
   &\leq\mathcal{M}_2 [\|u\|_{L_{p}^{\frac{n}{p}}\left(\mathbb{B}\right)}^{p}+1+\|u_t\|^2_{L^{\frac{n}{2}}_{2}\left(\mathbb{B}\right)}+\mathcal{H}(t)],
   \end{split}
\end{equation}
where
\begin{equation}\label{5.21}
  \mathcal{M}_2:=2^{\frac{\sigma}{1-\sigma}}\max\left\{1,\varepsilon^{\frac{\sigma}{1-\sigma}}C_1,
\varepsilon^{\frac{1}{1-\sigma}}C_2 {\frac{2}{p[2(1-\sigma)-1]}},\varepsilon^{\frac{1}{1-\sigma}}C_2 {\frac{p[2(1-\sigma)-1]-2}{p[2(1-\sigma)-1]}}\right\}.
\end{equation}\\
\textbf{Step 3: Blow UP.}\\
Take a combination \eqref{5.18}  with \eqref{5.20}, then $$ \mathcal{F}^ { \frac { 1 } { 1 - \sigma } } ( t ) \leq \frac { \mathcal{M}_2 } { \mathcal{M}_1 } \mathcal{F} ^ { \prime } ( t ) ,$$
 which implies by Gronwall's inequality
$$ \mathcal{F} ^ { \frac { \sigma } { 1 - \sigma } } ( t ) \geq \frac { 1 } { \mathcal{F} ^ { - \frac { \sigma } { 1 - \sigma } } ( 0 ) - \frac {\mathcal{M}_2 } {\mathcal{M}_1 } \frac { \sigma } { 1 - \sigma } t } ,$$
which yields $F(t)\rightarrow+\infty$ in finite time $T_{max}$ and
$$ T_{max} \leq \mathcal{F} ^ { - \frac { \sigma } { 1 - \sigma } } ( 0 ) \frac { \mathcal{M}_1 } { \mathcal{M}_2 } \frac { 1 - \sigma } { \sigma } .$$
where we also use the fact
\begin{equation}\label{5.22}
  \mathcal{F}(0)=\varepsilon \int_{\mathbb{B}}u_0u_1 \frac{dx_{1}}{x_{1}}dx'>0.
\end{equation}
by the assumption. One remains to prove that
$$\lim _ { t \rightarrow T_{max}^{-}  } \mathcal{F}(t)= + \infty\Longrightarrow \lim _ { t \rightarrow T_{max}^{-}  } \|u\|_{L_{p}^{\frac{n}{p}}\left(\mathbb{B}\right)}= + \infty.$$
\textbf{Case 1:} $H(t)\rightarrow+ \infty$\\
In this case, we have
$$\mathcal{H}(t)=E(0)-E(t)\leq E(0)-\frac{1}{p}\|g(x)^\frac{1}{p}u\|_{L_{p}^{\frac{n}{p}}\left(\mathbb{B}\right)}^{p}\leq
E(0)-\frac{\alpha}{p}\|u\|_{L_{p}^{\frac{n}{p}}\left(\mathbb{B}\right)}^{p},$$  which implies
$$\lim_ { t \rightarrow T^-_{max}  } \|u\|_{L_{p}^{\frac{n}{p}}\left(\mathbb{B}\right)}= + \infty.$$
\textbf{Case 2:} $\int_{\mathbb{B}}u u_t \frac{dx_{1}}{x_{1}}dx'\rightarrow+ \infty$\\
We have $$\int_{\mathbb{B}}u u_t \frac{dx_{1}}{x_{1}}dx'\leq\frac{1}{2}\|u_t\|^2_{L^{\frac{n}{2}}_{2}\left(\mathbb{B}\right)}
+\frac{1}{2}\|u\|^2_{L^{\frac{n}{2}}_{2}\left(\mathbb{B}\right)}
\leq\frac{1}{2}\|u_t\|^2_{L^{\frac{n}{2}}_{2}\left(\mathbb{B}\right)}
+\frac{1}{2\lambda_1}\|\nabla_{\mathbb{B}}u\|^{2}_{L_{2}^{\frac{n}{2}}\left(\mathbb{B}\right)}.$$
Recalling the $E(0)\geq E(t)$ and $E'(t)=-\|u_{t}\|^{m}_{L^{\frac{n}{m}}_{m}\left(\mathbb{B}\right)}$, we have
$$\frac{1}{2}\|u_t\|^2_{L^{\frac{n}{2}}_{2}\left(\mathbb{B}\right)}
+(1-\gamma(C^\ast)^{2})\|\nabla_{\mathbb{B}}u\|^{2}_{L_{2}^{\frac{n}{2}}\left(\mathbb{B}\right)}
\leq E(t)+\frac{1}{p}\|g(x)^\frac{1}{p}u\|_{L_{p}^{\frac{n}{p}}\left(\mathbb{B}\right)}^{p}
\leq E(0)+\frac{\beta}{p}\|u\|_{L_{p}^{\frac{n}{p}}\left(\mathbb{B}\right)}^{p},$$
which implies
$$\lim _ { t \rightarrow T^-_{max}  } \|u\|_{L_{p}^{\frac{n}{p}}\left(\mathbb{B}\right)}= + \infty.$$
This proof is complete.
\end{proof}

Finally,  we give a sufficient and necessary condition for the blowup result when the initial energy $E(0)<d$.

\begin{proposition}\label{Z}
Let all the assumptions in Theorem \ref{th5.2} hold. In addition, we assume $2<p<2+\frac{4}{n}$. Let $u(t)$ be a local solution to problem \eqref{1.1} on $[0,T_{max})$, if  $T_{max}<+\infty$, then $$\lim\limits_{t\to T_{max}^{-}}
\|u\|_{L^{\frac{n}{2}}_{2}\left(\mathbb{B}\right)}=+\infty.$$
\end{proposition}
\begin{proof}
Suppose that there exists an increasing sequence $\left\{t_{m}\right\}$ with $t_m\rightarrow T_{max}$
(as $m\to +\infty$) such that
$$\|u(t_m)\|_{L^{\frac{n}{2}}_{2}\left(\mathbb{B}\right)}\leq C $$ for some constants $C>0$.
It follows from \eqref{2.3}, \eqref{2.5} and \eqref{2.7} that
\begin{equation} \label{5.23}
\|\nabla_{\mathbb{B}}u\|^{2}_{L_{2}^{\frac{n}{2}}\left(\mathbb{B}\right)}
-\gamma\|V(x)^{\frac{1}{2}}u\|^{2}_{L_{2}^{\frac{n}{2}}\left(\mathbb{B}\right)}
\leq 2E(0)
+\frac{2}{p}\|g(x)^\frac{1}{p}u\|_{L_{p}^{\frac{n}{p}}\left(\mathbb{B}\right)}^{p},
\end{equation}
for all $t\in[0,T_{max})$.
Next, from Lemma \ref{GN} and Young's inequality, we have
\begin{equation*}
\begin{split}
\|g(x)^\frac{1}{p}u\|_{L_{p}^{\frac{n}{p}}\left(\mathbb{B}\right)}^{p}
   &\leq C\beta\|u\|^{p(1-\theta)}
   _{L^{\frac{n}{2}}_{2}\left(\mathbb{B}\right)}
   \|\nabla_{\mathbb{B}}u\|^{p\theta}
   _{L_{2}^{\frac{n}{2}}\left(\mathbb{B}\right)}\\[2mm]
   &\leq \beta C \frac{\varepsilon^{q_1}}{q_1}
   \|\nabla_{\mathbb{B}}u\|^{p\theta q_1}
   _{L_{2}^{\frac{n}{2}}\left(\mathbb{B}\right)}
   +\beta C\frac{1}{q_2\varepsilon^{q_2}}\|u\|^{p(1-\theta)q_2}
   _{L^{\frac{n}{2}}_{2}\left(\mathbb{B}\right)}
\end{split}
\end{equation*}
with some positive constant $C$, where $\theta=\frac{(p-2)n}{2p}$ and $\varepsilon>0$ is a constant to be determined, $q_1, q_2>0$ and $\frac{1}{q_1}+\frac{1}{q_2}=1$ . Thus, from \eqref{5.23} we have
\begin{equation*}
\begin{split}
&\|\nabla_{\mathbb{B}}u(t)\|^{2}_{L_{2}^{\frac{n}{2}}\left(\mathbb{B}\right)}
-\gamma\|V(x)^{\frac{1}{2}}u(t)\|^{2}_{L_{2}^{\frac{n}{2}}\left(\mathbb{B}\right)}\\[2mm]
&\leq 2E(0)+\frac{2\beta C}{p}[\frac{\varepsilon^{q_1}}{q_1}
   \|\nabla_{\mathbb{B}}u(t)\|^{p\theta q_1}
   _{L_{2}^{\frac{n}{2}}\left(\mathbb{B}\right)}
   +\frac{1}{q_2\varepsilon^{q_2}}\|u(t)\|^{p(1-\theta) q_2}
   _{L^{\frac{n}{2}}_{2}\left(\mathbb{B}\right)}].
 \end{split}
\end{equation*}
Taking $t=t_m$ and $q_1=\frac{4}{\left(p-2\right)n}$, we get
\begin{equation*}
\begin{split}&\|\nabla_{\mathbb{B}}u\left(t_m\right)\|^{2}_{L_{2}^{\frac{n}{2}}\left(\mathbb{B}\right)}
-\gamma\|V(x)^{\frac{1}{2}}u\left(t_m\right)\|^{2}_{L_{2}^{\frac{n}{2}}\left(\mathbb{B}\right)}
\leq C_1 +\frac{2\beta C\varepsilon^{q_1}}{q_1}
\|\nabla_{\mathbb{B}}u\left(t_m\right)\|^{2}_{L_{2}^{\frac{n}{2}}\left(\mathbb{B}\right)}
\end{split}
\end{equation*}
Choosing $\varepsilon$ sufficiently  small, by Lemma \ref{d}, we have
$$\|\nabla_{\mathbb{B}}u\left(t_m\right)\|_{L_{2}^{\frac{n}{2}}\left(\mathbb{B}\right)}
\leq\frac{{C_2}^2}{1-\gamma(C^\ast)^{2}}.$$
Combining  Lemma \ref{b},  we have
$$\|g(x)^\frac{1}{p}u\left(t_m\right)\|_{L_{p}^{\frac{n}{p}}\left(\mathbb{B}\right)}
\leq C_*\|\nabla_{\mathbb{B}}u\left(t_m\right)\|^{2}_{L_{2}^{\frac{n}{2}}\left(\mathbb{B}\right)}\leq C.$$
Using Lemma \ref{e}, we have $\|u_t\left(t_m\right)\|^2_{L^{\frac{n}{2}}_{2}\left(\mathbb{B}\right)}$
is bounded. Thus,
$$\lim_{m\rightarrow+\infty}[\|u_t\left(t_m\right)\|^2_{L^{\frac{n}{2}}_{2}\left(\mathbb{B}\right)}
+\|\nabla_{\mathbb{B}}u\left(t_m\right)\|^{2}_{L_{2}^{\frac{n}{2}}\left(\mathbb{B}\right)}]
<\infty$$
which  contradicts with $T_{max}<+\infty$.
\end{proof}

\begin{theorem}\label{th5.5}
Let the conditions in Proposition \ref{Z} hold,  and $u(t)$ be a local solution to \eqref{1.1} on $[0,T_{max})$.    Then  there must exist a real number $t_0\in[0,T_{max})$ such that $E(t_0)<d$ and $u(t_0)\in\mathcal{V}$ if and only if $T_{max}<\infty$, i.e. the solution blows up at a finite time.
\begin{proof}
The `` if part" is a direct result of Theorem \ref{th5.2}. So we only need to prove just the ``only if part". In fact, it is easy to see that
$$|\|u(t)\|_{L^{\frac{n}{2}}_{2}\left(\mathbb{B}\right)}-\|u_0\|_{L^{\frac{n}{2}}_{2}\left(\mathbb{B}\right)}|
\leq\|u(t)-u_0\|_{L^{\frac{n}{2}}_{2}\left(\mathbb{B}\right)}
\leq\int^t_0\|u_{\tau}(\tau)\|_{L^{\frac{n}{2}}_{2}\left(\mathbb{B}\right)}d\tau$$
on $[0,T_{max})$.
Using  H\"{o}lder inequality, we obtain
$$|\|u(t)\|_{L^{\frac{n}{2}}_{2}\left(\mathbb{B}\right)}-\|u_0\|_{L^{\frac{n}{2}}_{2}\left(\mathbb{B}\right)}|
\leq C\int^t_0\|u_t(\tau)\|_{L^{\frac{n}{m}}_{m}\left(\mathbb{B}\right)}d\tau
\leq C t^{\frac{m-1}{m}}\left(\int^t_0\|u_t(\tau)\|^m_{L^{\frac{n}{m}}_{m}\left(\mathbb{B}\right)}\right)^{\frac{1}{m}}.$$
Consequently, it follows from Lemma \ref{e} that
\begin{equation}\label{5.24}
E(0)-E(t)=\int^t_0\|u_{\tau}(\tau)\|^m_{L^{\frac{n}{m}}_{m}\left(\mathbb{B}\right)}d\tau
\geq Ct^{1-m}|\|u(t)\|_{L^{\frac{n}{2}}_{2}\left(\mathbb{B}\right)}
-\|u_0\|_{L^{\frac{n}{2}}_{2}\left(\mathbb{B}\right)}|^m.
\end{equation}
Then it follows from Proposition \ref{Z} that $\lim\limits_{t\to T_{max}^{-}}\|u(t)\|_{L^{\frac{n}{2}}_{2}\left(\mathbb{B}\right)}=+\infty$, so \eqref{5.24} implies
\begin{equation}\label{5.25}
\lim\limits_{t\rightarrow T_{max}^-}E(t)=-\infty.
\end{equation}
On the other hand, noticing
\begin{equation*}
E(0)\geq E(t)\geq\frac{1}{2}\|u_t\|^2_{L^{\frac{n}{2}}
_{2}\left(\mathbb{B}\right)}+\frac{1}{p}I(u)
  +\frac{(p-2)}{2p}(1-\gamma(C^\ast)^{2})
  \|\nabla_{\mathbb{B}}u\|^{2}_{L_{2}^{\frac{n}{2}}\left(\mathbb{B}\right)}.
\end{equation*}
and  $\lim\limits_{t \to T^-_{max}}[\|u_t\left(t\right)\|^2_{L^{\frac{n}{2}}_{2}\left(\mathbb{B}\right)}
+\|\nabla_{\mathbb{B}}u\left(t\right)\|^{2}_{L_{2}^{\frac{n}{2}}\left(\mathbb{B}\right)}]
<\infty$, we have
\begin{equation}\label{5.26}
\lim_{t\rightarrow T_{max}^-}I(u(t))=-\infty.
\end{equation}
By \eqref{5.25} and \eqref{5.26}, we obtain that
$$J(u\left(t_0\right))\leq E(t_0)<d$$
and
$$I(u\left(t_0\right))<0$$
for some $t_0\in[0,T_{max})$, which implies that $u\left(t_0\right)\in\mathcal{V}$
and $E(t_0)<d$. The proof is complete.
\end{proof}
\end{theorem}

\section*{Acknowledgments}
\noindent The authors would like to thank the referees for the careful reading of this paper and for the valuable suggestions to
improve the presentation and the style of the paper. The project is supported by Natural Science Foundation of Henan
Province (No.242300420241),  and the Open project of
Complexity Science Center of Henan University of Technology (No. CSKFJJ-2024-9).

\end{document}